# LOCAL ALIGNMENT OF MARKOV CHAINS

BY NIELS RICHARD HANSEN

*University of Copenhagen*

We consider local alignments without gaps of two independent Markov chains from a finite alphabet, and we derive sufficient conditions for the number of essentially different local alignments with a score exceeding a high threshold to be asymptotically Poisson distributed. From the Poisson approximation a Gumbel approximation of the maximal local alignment score is obtained. The results extend those obtained by Dembo, Karlin and Zeitouni [*Ann. Probab.* **22** (1994) 2022–2039] for independent sequences of i.i.d. variables.

**1. Introduction.** Local alignment of two biological sequences (DNA-molecules or proteins) is one of the most important and used tools in modern molecular biology for locating highly similar contiguous parts of the sequences. High similarity is usually interpreted as an evolutionary or functional relationship between the molecules. We show how the distribution of local alignment similarity scores behaves asymptotically when aligning independent Markov chains.

It is important to understand the distribution of local alignment scores for assessing the significance of, for example, the maximally scoring local alignment. Formally this is a test of the null hypothesis that two sequences are independent Markov chains against a somewhat unspecified alternative that they are not independent. The test statistic considered is the maximal local similarity score.

Usually when considering local alignments we are interested in not only the maximally scoring local alignment but also other essentially different local alignments that reach a score above a given threshold. It is therefore useful also to know the distribution of the number of local alignments of independent Markov chains that reach a score above a given threshold. In fact, it is this problem that we handle in the first place and the obtained

Poisson approximation can easily be turned into a Gumbel approximation of the distribution of the maximal local alignment score.

The kind of local alignment we consider is *gapless* local alignment meaning that we search for (contiguous) parts of the two sequences that attain a high similarity when matched letter by letter. Similarity is measured by adding up a score for each pair of matched letters. In practice it is common to allow for the insertion of gaps in the sequences—each gap adding a suitable penalty to the similarity score—which usually increases the power of the test. The introduction of gaps does, however, make the problem of understanding the asymptotic distribution of local alignment scores substantially more complicated although progress for i.i.d. sequences has been made more recently; see [4, 10, 17, 18]. In another direction, exact distributional results for i.i.d. sequences can be obtained if "shifting" is not allowed and if the scores are integer valued; see [14]. This work has also been generalized to Markov sequences; see [15].

The main result is stated as Theorem 3.1. It says that if the expected similarity score under the null hypothesis is negative, then there exist constants $\theta^*, K^* > 0$ such that if we let $s$ denote the maximal local alignment score obtained when aligning two independent Markov chains of length $n$, then the normalized score defined by

$$s' = \theta^* s - \log(K^* n^2) \tag{1}$$

approximately follows a Gumbel distribution for $n \to \infty$. Moreover, the number of normalized local alignment scores exceeding the threshold $x$ is approximately Poisson distributed with mean $\exp(-x)$ for $n \to \infty$. We have ignored some details and there are certain assumptions that need to be fulfilled for this to be a mathematically rigorous statement. We refer to Theorem 3.1 and its prerequisites.

It should be mentioned that the results are the expected generalizations of those obtained in [6] for independent i.i.d. sequences, but the techniques of proof are not straightforward generalizations. Indeed, this author would like to emphasize the novelty of certain techniques developed in this paper. In particular the results achieved in Sections 5.4 and 5.5 may be of independent interest. Moreover, the framework of Markov chains does not only provide a change of the null hypothesis but it also opens up the possibility of choosing new types of score functions as we discuss in Remark 3.5. This can increase the power of the test. In addition, by expanding the state space suitably the results obtained in this paper also cover null hypotheses where the aligned sequences have a higher-order Markov dependency or come from a hidden Markov model.

**2. Local gapless alignment.** Let $(X_n)_{n\geq 1}$ and $(Y_n)_{n\geq 1}$ be two sequences of random variables taking values in a finite set $E$. We compare parts of one sequence with parts of the other using a score function $f: E \times E \to \mathbb{Z}$, and we define the random variables

$$S_{i,j}^\delta = \sum_{k=1}^{\delta} f(X_{i+k}, Y_{j+k}),$$

for $i, j, \delta \geq 0$. The variable $S_{i,j}^\delta$ is the local score for the local comparison of the sequence part $X_{i+1} \cdots X_{i+\delta}$ with the sequence part $Y_{j+1} \cdots Y_{j+\delta}$.

We make the assumption that $f$ takes integer values to emphasize the lattice nature of $f$ that is often met in practice. To assure that $\mathbb{Z}$ indeed is the minimal lattice, the greatest common divisor of the integers $f(x, y)$, $x, y \in E$, is assumed to be 1. The results obtained are valid if $f$ takes real, nonlattice values in a slightly modified form; see Remark 3.2.

The score function can be regarded as an $E \times E$ matrix, which is convenient when writing down the values $f(x, y)$. We will find it most useful to simply regard $f$ as an element in a vector space. Probability measures will then be regarded as elements in the dual space and we use the functional notation

$$\nu(f) = \sum_{x,y} f(x,y) \nu(x,y)$$

to denote the mean of $f$ evaluated under the probability measure $\nu$.

For $n \geq 1$ define

$$\mathcal{H}_n = \{(i, j, \delta) | 0 \leq i \leq i + \delta \leq n, \, 0 \leq j \leq j + \delta \leq n\}.$$

The elements $(i, j, \delta) \in \mathcal{H}_n$ are called alignments.

We want to understand the distribution of the collection

$$(S_{i,j}^\delta)_{(i,j,\delta) \in \mathcal{H}_n}$$

of local scores over all alignments. We will in particular be interested in the distribution of

(2) $$\mathcal{M}_n = \max_{(i,j,\delta) \in \mathcal{H}_n} S_{i,j}^\delta,$$

the maximal local score over the set of alignments. We will also study the number, $C_n(t)$, say, of *essentially different* variables $S_{i,j}^\delta$ in $\mathcal{H}_n$ exceeding some threshold $t \geq 0$. What we mean by "essentially different" is defined precisely below.

The local scores are efficiently summarized in the score matrix $(T_{i,j})_{0 \leq i,j \leq n}$, which is defined as follows. For $i = 0$ or $j = 0$ let $T_{i,j} = 0$ and define recursively

(3) $$T_{i,j} = \max\{T_{i-1,j-1} + f(X_i, Y_j), 0\}$$

for $i,j \geq 1$. As we will show (cf. Remark 3.6 below), the maximum $\mathcal{M}_n$ can be computed as

$$\mathcal{M}_n = \max_{i,j} T_{i,j}. \tag{4}$$

This fact is closely related to the idea in the Smith–Waterman algorithm for computing the (gapped) maximal local alignment score efficiently; see [21].

DEFINITION 2.1. For $0 \leq i,j \leq n-1$ define

$$\Delta(i,j) = \inf\{\delta > 0 | S_{i,j}^\delta \leq 0, \text{ or } i+\delta = n, \text{ or } j+\delta = n\}.$$

If $T_{i,j} = 0$, the alignment $(i,j,\Delta(i,j))$ is called an excursion, and we let $\mathcal{E}_n$ denote the set of all excursions.

Note that $\mathcal{E}_n$ is a stochastic subset of $\mathcal{H}_n$. It follows from the definition of the score matrix $(T_{i,j})$ and the definition of an excursion that if $(i,j,\Delta) \in \mathcal{E}_n$ and $0 < \delta < \Delta$, then

$$T_{i+\delta,j+\delta} = S_{i,j}^\delta.$$

An excursion corresponds to a diagonal strip in the score matrix, for which the score starts at zero and then stays strictly positive along that diagonal strip until it either reaches zero or the indices hit the boundary of the score matrix.

The maximum over an excursion $e = (i,j,\Delta) \in \mathcal{E}_n$ is denoted by

$$\mathcal{M}_e = \max_{0 < \delta \leq \Delta} T_{i+\delta,j+\delta}. \tag{5}$$

DEFINITION 2.2. The number of essentially different excesses over $t$ is defined as

$$C_n(t) = \sum_{e \in \mathcal{E}_n} 1(\mathcal{M}_e > t). \tag{6}$$

From (4) it follows that $(C_n(t) = 0) = (\mathcal{M}_n \leq t)$.

**3. Alignment of independent Markov chains.** Assume that the stochastic processes $(X_n)_{n \geq 1}$ and $(Y_n)_{n \geq 1}$ are independent Markov chains with transition probabilities $P$ and $Q$, respectively. Assume that $P$ and $Q$ are irreducible and aperiodic matrices with left invariant probability vectors $\pi_P$ and $\pi_Q$, respectively. Let $\pi = \pi_P \otimes \pi_Q$. With

$$\mu = \pi(f) = \sum_{x,y \in E} f(x,y) \pi_P(x) \pi_Q(y)$$

the (invariant) mean of $f(X_1,Y_1)$ we will assume throughout that $\mu < 0$.

In the following, a *cycle* w.r.t. a matrix of transition probabilities $P$ is a finite sequence $x_1, \ldots, x_n$ such that

$$P(x_i, x_{i+1 (\mathrm{mod}\, n)}) > 0$$

for $i = 1, \ldots, n$. We will assume that the following regularity conditions on $f$, $P$ and $Q$ are fulfilled: For some $n \geq 1$ there exist cycles $x_1, \ldots, x_n$ (w.r.t. $P$) and $y_1, \ldots, y_n$ (w.r.t. $Q$) such that

(7) $$\sum_{k=1}^{n} f(x_k, y_k) > 0.$$

For any $T \geq 1$ there exist an $n \geq 1$ and cycles $x_1, \ldots, x_n$ (w.r.t. $P$) and $y_1, \ldots, y_n$ (w.r.t. $Q$) such that

(8) $$\sum_{k=1}^{n} f(x_k, y_k) \neq \sum_{k=1}^{n} f(x_k, y_{k+T (\mathrm{mod}\, n)}).$$

See Remark 3.3 below for comments related to this somewhat strange looking condition.

For convenience we will assume that both Markov chains are stationary, though the results obtained hold anyway. We denote by $\mathbb{P}$ the probability measure $\mathbb{P}_\pi$ under which $(X_n, Y_n)_{n \geq 1}$ is a stationary Markov chain with transition probabilities $P \otimes Q$. It will in addition be convenient to assume that there exist auxiliary random variables $X_0$ and $Y_0$ such that $(X_n, Y_n)_{n \geq 0}$ under $\mathbb{P}$ forms a stationary Markov chain too. As usual $\mathbb{P}_{x,y}$ will denote the probability measure where $X_1 = x$ and $Y_1 = y$.

We define for $\theta \in \mathbb{R}$ an $E^2 \times E^2$ matrix $\Phi(\theta)$ with positive entries by

$$\Phi(\theta)_{(x,y),(x',y')} = \exp(\theta f(x', y')) P_{x,x'} Q_{y,y'},$$

and we let $\varphi(\theta)$ denote the spectral radius (the Perron–Frobenius eigenvalue) of this matrix. Then $\varphi$ is a convex $C^\infty$-function in $\theta$, and due to (7), $\varphi(\theta) \to \infty$ for $\theta \to \infty$. The fact that $\varphi$ is (log)convex is due to Kingman [12], and the implicit function theorem can be used to show that $\varphi$ is $C^\infty$. Furthermore, by Corollary XI.2.9(a) in [3] it holds that

(9) $$\partial_\theta \varphi(0) = \mu,$$

hence if $\mu < 0$ there exists a (by convexity unique) solution $\theta^* > 0$ to the equation $\varphi(\theta) = 1$. If $r^*$ denotes the (up to scaling unique) right eigenvector corresponding to the eigenvalue 1 for $\Phi(\theta^*)$, the matrix defined by

$$R^*_{(x,y),(x',y')} = \frac{r^*(x', y')}{r^*(x, y)} \Phi(\theta^*)_{(x,y),(x',y')}$$

is an irreducible stochastic matrix with a unique left invariant probability vector, which we will denote by $\pi^*$.

With $g: E^2 \times E^2 \to \mathbb{R}$ any given function we introduce two $E^3 \times E^3$ matrices, $\Phi_1(g)$ and $\Phi_2(g)$, by

$$\Phi_1(g)_{(x,y,z),(x',y',z')} = \exp(g(x,y,x',y') + g(x,z,x',z'))P_{x,x'}Q_{y,y'}Q_{z,z'},$$

$$\Phi_2(g)_{(x,w,y),(x',w',y')} = \exp(g(x,y,x',y') + g(w,y,w',y'))P_{x,x'}P_{w,w'}Q_{y,y'},$$

and we let $\varphi_1(g)$ and $\varphi_2(g)$ denote the corresponding spectral radii. In terms of the functions $\varphi_1$ and $\varphi_2$ we define

$$J_i = \sup_g \{2\hat{\pi}(g) - \log \varphi_i(g)\}$$

for $i = 1, 2$. Here $\hat{\pi} = \pi^* \otimes R^*$ denotes the measure on $E^2 \times E^2$ with point probabilities $\hat{\pi}(x, y, x', y') = \pi^*(x, y) R^*_{(x,y),(x',y')}$. We discuss $J_1$ and $J_2$ in further detail in Remark 3.8.

Finally, if we define the process $(S_n)_{n \geq 0}$ by $S_0 = 0$ and for $n \geq 1$,

$$(10) \qquad S_n = \sum_{k=1}^{n} f(X_k, Y_k),$$

we can define a constant, $K^*$, in terms of this process, as done, for example, by (1.26) in Theorem B in [11]. We discuss this constant in further detail in Remark 3.7.

THEOREM 3.1. *Assume that $\mu < 0$, that the regularity conditions given by (7) and (8) are fulfilled, and that $\theta^*$ and $K^*$ are the constants defined above. Define for $x \in \mathbb{R}$*

$$(11) \qquad t_n = \frac{\log K^* + \log n^2 + x}{\theta^*}$$

*and $x_n \in [0, \theta^*)$ by $x_n = \theta^*(t_n - \lfloor t_n \rfloor)$. Then if*

$$(12) \qquad 2\min\{J_1, J_2\} > 3\theta^* \pi^*(f),$$

*it holds that*

$$(13) \qquad \|\mathcal{D}(C_n(t_n)) - \mathrm{Poi}(\exp(-x + x_n))\| \to 0$$

*for $n \to \infty$. Here $\|\cdot\|$ denotes the total variation norm and $\mathcal{D}(C_n(t_n))$ is the distribution of $C_n(t_n)$. In particular*

$$(14) \qquad \mathbb{P}(\mathcal{M}_n \leq t_n) - \exp(-\exp(-x + x_n)) \to 0$$

*for $n \to \infty$.*

The theorem deserves a number of remarks.

REMARK 3.2. The choice of $x_n = \theta^*(t_n - \lfloor t_n \rfloor)$ assures that $t_n - x_n/\theta^* = \lfloor t_n \rfloor \in \mathbb{Z}$. Due to the lattice effect arising from $f$ taking values in $\mathbb{Z}$ it follows that

$$(C_n(t_n) = m) = (C_n(t_n - x_n/\theta^*) = m)$$

as well as

$$(\mathcal{M}_n \leq t_n) = (\mathcal{M}_n \leq t_n - x_n/\theta^*),$$

and this is the reason that we need to correct by $x_n$ in the asymptotic formulas. If $f$ is a real, nonlattice function, Theorem 3.1 holds without the $x_n$-correction.

REMARK 3.3. The regularity condition (8) does not look particularly nice in general but is usually satisfied by quite trivial arguments. Essentially we want to avoid the situation where

(15) $$f(x,y) = f_1(x) + f_2(y)$$

for two functions $f_1, f_2 : E \to \mathbb{R}$. It is clear that if $f$ is of the form (15), then (8) does *not* hold. It is easy to verify that if $P$ and $Q$ have only strictly positive entries, condition (8) is equivalent to $f$ *not* being of the form (15). In general, however, this author has not been able to prove that $f$ *not* being of the form (15) is sufficient for (8) to hold. On the other hand, no counterexamples have been found either. In the proof we will explicitly need that (8) holds.

REMARK 3.4. It is possible, and of practical relevance, to allow for the aligned sequences to have different lengths $m$ and $n$, say. In this case Theorem 3.1 holds for $n, m \to \infty$ with

$$t_{m,n} = \frac{\log K^* + \log(mn) + x}{\theta^*}.$$

Some restriction on the simultaneous growth of $m$ and $n$ must be made in order for this to be true. In the proof of Lemma 5.15 we will need to be able to choose integers $l_{n,m}$ fulfilling that

$$\lim_{n,m \to \infty} \frac{\log(nm)}{l_{n,m}} = \lim_{n,m \to \infty} \frac{l_{n,m}}{\min\{n,m\}} = 0,$$

where $n, m \to \infty$ refers to the desired simultaneous growth of $n$ and $m$. Clearly this can be achieved if $m \sim cn$ for some constant $c > 0$, whereas, for example, $m \sim \log n$ does not work.

REMARK 3.5. For notational convenience Theorem 3.1 was stated and proved using a score function $f$ that depends on a single pair of variables only. When aligning Markov chains it would be perfectly natural to use a score function that depends on *pair-transitions* instead, that is, $f : E^2 \times E^2 \to \mathbb{R}$ and

$$S_{i,j}^\delta = \sum_{k=1}^\delta f(X_{i+k-1}, Y_{j+k-1}, X_{i+k}, Y_{j+k}).$$

Theorem 3.1 holds for this kind of score function with the obvious modifications. For instance, $\Phi$ is defined as

$$\Phi(\theta)_{(x,y),(x',y')} = \exp(\theta f(x,y,x',y')) P_{x,x'} Q_{y,y'},$$

and $\pi^*$ in (12) is replaced by $\hat{\pi}$. In practice $f$ can be chosen as a (conditional) log-likelihood ratio. If the *alternative* to the null hypothesis is assumed to be a Markov chain on $E^2$ governed by an $E^2 \times E^2$ matrix of transition probabilities $R$, then we could choose

$$f(x,y,x',y') = \log \frac{R_{(x,y),(x',y')}}{P_{x,x'} Q_{y,y'}}.$$

This score function does clearly not take integer values in general, but one may choose to consider $\lfloor Nf \rfloor$ for suitably large $N$ if integer scores are preferred.

We find that for this score function $f$ and for $\theta = 1$

$$\begin{aligned}\Phi(1)_{(x,y),(x',y')} &= \exp(f(x,y,x',y')) P_{x,x'} Q_{y,y'} \\ &= R_{(x,y),(x',y')},\end{aligned}$$

which has row sums equal to 1. Hence $\varphi(1) = 1$ implying that $\theta^* = 1$.

REMARK 3.6. The process $(S_n)_{n \geq 0}$ defined by (10) is called a Markov controlled random walk or a Markov additive process (abbreviated MAP); see [3], Chapter XI. The process $(T_n)_{n \geq 0}$ defined by

(16) $$T_n = S_n - \min_{0 \leq k \leq n} S_k$$

is called the reflection of the MAP at the zero barrier. It is straightforward to verify that $(T_n)_{n \geq 0}$ satisfies the recursion

$$T_n = \max\{T_{n-1} + f(X_n, Y_n), 0\}$$

for $n \geq 1$. In addition

(17) $$\begin{aligned}\max_{1 \leq k \leq m \leq n} S_m - S_k &= \max_{1 \leq m \leq n} \left\{ S_m - \min_{1 \leq k \leq m} S_k \right\} \\ &= \max_{1 \leq m \leq n} T_m.\end{aligned}$$

We see that $S_{0,0}^n = S_n$ and $T_{n,n} = T_n$. Thus along the main diagonal in the score matrix $(T_{i,j})_{0 \leq i,j \leq n}$ we find the reflection of the MAP $(S_n)_{n \geq 0}$. Along all other diagonals in the score matrix we find the reflections of MAPs too—these MAPs being defined by shifting the Markov chain $(X_n)_{n \geq 1}$ along $(Y_n)_{n \geq 1}$. It follows from (17) that (4) indeed holds. Due to independence and stationarity of the two Markov chains all the reflected MAPs along diagonals have the same distribution, but they are dependent. The interpretation of Theorem 3.1 is that asymptotically the number of excursions exceeding level $t_n$ has the same distribution as if the reflected MAPs were independent.

REMARK 3.7. The constant $K^*$ is defined in terms of the MAP $(S_n)_{n \geq 0}$. Let $\tau_-(0) = 0$ and for $k \geq 1$ let

$$\tau_-(k) = \inf\{n > \tau_-(k-1) | S_n \leq S_{\tau_-(k-1)}\}$$

denote the times when the MAP descends below its previous minimum. These stopping times are known as the descending ladder epochs for the MAP, and they are almost surely finite due to assumption that $\mu < 0$. One should note that $\tau_-(k)$ is also the $k$th time that the reflected MAP $(T_n)_{n \geq 0}$ hits 0. A thorough treatment of the ladder epochs is given in [1] covering also general state-space Markov chains. From Theorem 1(i) in [1] it follows that the sampled Markov chain $(X_{\tau_-(n)}, Y_{\tau_-(n)})_{n \geq 0}$ has a unique invariant probability distribution, which we will denote by $\nu$. As we consider only a finite state-space Markov chain, this is also a direct consequence of the Wiener–Hopf factorization ([3], Theorem XI.2.12). Moreover, the sequence defined by

$$u_{x,y}(n) = \mathbb{P}(T_n = 0, X_n = x, Y_n = y)$$
$$= \mathbb{P}(\exists k : \tau_-(k) = n, X_n = x, Y_n = y)$$

for $x, y \in E$ and $n \geq 1$ forms a renewal sequence and the elementary renewal theorem, ([3], Theorem V.1.4) gives that

$$(18) \qquad \frac{1}{n} \sum_{k=1}^{n} u_{x,y}(k) \to \frac{\nu(x,y)}{\mu_-}$$

for $n \to \infty$ where $\mu_- = \mathbb{E}_\nu(\tau_-(1))$. We refer to ([5] Theorem 10.4.3) for a proof that the inverse of the mean recurrence time indeed is given as the right-hand side limit above.

As stated in Lemma B in [11], when $\mu < 0$ and (7) holds, then

$$(19) \qquad \lim_{u \to \infty} \mathbb{P}_{x,y}\left(\max_{1 \leq n \leq \tau_-(1)} S_n > u\right) \exp(\theta^* u) = e(x,y)$$

for some constants $e(x,y) \geq 0$, $x, y \in E$. In terms of these limits the constant $K^*$ can be represented as

$$K^* = \frac{1}{\mu_-} \sum_{x,y} \nu(x,y) e(x,y). \tag{20}$$

As a consequence of Walds identity for MAPs ([3], Corollary XI.2.6), it holds that

$$\mu_- = \frac{\mathbb{E}_\nu(S_{\tau_-(1)})}{\mu},$$

which shows that (20) is identical to the representation of $K^*$ in (1.26) in [11]. We refer to [11] for more details and in particular their Section 5 for issues related to the computation of $K^*$.

REMARK 3.8. The function $g \mapsto \log \varphi_i(g)$ is a convex function and $J_i$ is thus the Fenchel–Legendre transform of the function evaluated in $2\hat{\pi}$. It is possible to identify $J_i$ as the value of a large deviation rate-function. Considering $J_1$ we introduce the function $h \colon E^3 \times E^3 \to \mathbb{R}^{E^2 \times E^2}$ by

$$h(x,y,z,x',y',z') = (\mathbf{1}_{(x,y),(x',y')}(v) + \mathbf{1}_{(x,z),(x,z')}(v))_{v \in E^2 \times E^2}.$$

If $(X_n, Y_n, Z_n)_{n \geq 0}$ is a Markov chain with transition probabilities $P \otimes Q \otimes Q$, the large deviation rate-function for the empirical average

$$\frac{1}{n} \sum_{k=1}^{n} h(X_{k-1}, Y_{k-1}, Z_{k-1}, X_k, Y_k, Z_k)$$

evaluated in $2\alpha$ for $\alpha$ a probability measure on $E^2 \times E^2$ is given as

$$I(2\alpha) = \sup_g \{2\alpha(g) - \log \varphi_1(g)\};$$

see Theorem 3.1.2 in [7]. In particular $J_1 = I(2\hat{\pi})$.

Let $\nu$ be a probability measure on the space $E^3 \times E^3$ and define

$$\nu^1(x,y,z) = \sum_{x',y',z'} \nu(x,y,z,x',y',z'),$$

$$\nu^2(x',y',z') = \sum_{x,y,z} \nu(x,y,z,x',y',z');$$

thus $\nu^1$ and $\nu^2$ are marginal probability measures on $E^3$. The measure $\nu$ is called *shift-invariant* if $\nu^1 = \nu^2$, and we denote by $\widetilde{\mathcal{M}}$ the set of shift-invariant probability measures on $E^3 \times E^3$. Considering the Markov chain on $E^3$ with transition probability matrix $P \otimes Q \otimes Q$, then the large deviation rate-function for the pair-empirical measure (cf. Theorem VI.3 in [8]) is given as

$$I^2(\nu) = \sum_{\substack{x,y,z \\ x',y',z'}} \nu(x,y,z,x',y',z') \log \frac{\nu(x,y,z,x',y',z')}{\nu^1(x,y,z)P_{x,x'}Q_{y,y'}Q_{z,z'}}$$

$$= H(\nu|\nu^1 \otimes P \otimes Q \otimes Q)$$

for $\nu \in \widetilde{\mathcal{M}}$. Here $H(\cdot|\cdot)$ denotes the relative entropy. Defining additional marginals

$$\nu_{12}(x,y,x',y') = \sum_{z,z'} \nu(x,y,z,x',y',z'),$$

$$\nu_{13}(x,z,x',z') = \sum_{y,y'} \nu(x,y,z,x',y',z'),$$

with $\nu_{12}$ and $\nu_{13}$ being probability measures on $E^2 \times E^2$, it is a consequence of the contraction principle, Theorem III.20 in [8], that

$$J_1 = \inf\{I^2(\nu)|\nu \in \widetilde{\mathcal{M}} : \nu_{12} + \nu_{13} = 2\hat{\pi}\}.$$

A similar representation of $J_2$ is of course possible.

If $(X_n)_{n\geq 1}$ and $(Y_n)_{n\geq 1}$ are independent sequences of i.i.d. variables with the $X$'s having distribution $\pi_1$ and the $Y$'s having distribution $\pi_2$, then

$$\varphi(\theta) = \mathbb{E}(\exp(\theta f(X_1,Y_1)))$$

is the Laplace transform of the distribution of $f(X_1,Y_1)$, and $\theta^* > 0$ solves $\varphi(\theta) = 1$. Moreover, $\pi^*$ is the probability measure on $E \times E$ with point probabilities $\pi^*(x,y) = \exp(\theta^* f(x,y))\pi_1(x)\pi_2(y)$. In this case we can verify that the infimum above is attained for

$$\nu(x,y,z,x',y',z') = \frac{\pi^*(x,y)\pi^*(x,z)\pi^*(x',y')\pi^*(x',z')}{\pi_1^*(x)\pi_1^*(x')}$$

with $\pi_1^*$ denoting the first marginal of $\pi^*$. To see this first note that $\nu$ is clearly shift-invariant with the desired marginal property, $\nu_{12} + \nu_{13} = 2\pi^* \otimes \pi^*$. A simple computation reveals that

$$I^2(\nu) = 2\theta^*\pi^*(f) - \sum_x \pi_1^*(x) \log \frac{\pi_1^*(x)}{\pi_1(x)} = 2\theta^*\pi^*(f) - H(\pi_1^*|\pi_1),$$

and for any other shift-invariant $\tilde{\nu}$ with the same marginal property one finds that

$$I^2(\tilde{\nu}) = H(\tilde{\nu}|\nu) + I^2(\nu),$$

where $H(\tilde{\nu}|\nu) \geq 0$, hence

$$J_1 = 2\theta^*\pi^*(f) - H(\pi_1^*|\pi_1).$$

Similarly we can show that $J_2 = 2\theta^*\pi^*(f) - H(\pi_2^*|\pi_2)$. Since $\theta^*\pi^*(f) = H(\pi^*|\pi_1 \otimes \pi_2)$, the condition given by (12) is equivalent to

(21) $$H(\pi^*|\pi_1 \otimes \pi_2) > 2\max\{H(\pi_1^*|\pi_1), H(\pi_2^*|\pi_2)\},$$

which is precisely the condition (E′) required in [6] in the i.i.d. case for $C_n(t_n)$ to be asymptotically Poisson distributed. Since condition (H) in [6] is equivalent to $\mu < 0$ and (7) in the i.i.d. case, and since condition (E′) actually implies that $f$ does *not* take the form (15) in the i.i.d. setup, we conclude that (8) is also fulfilled in the i.i.d. case when assuming (E′); see Remark 3.3. Thus Theorem 3.1 specializes in the i.i.d. case to Theorem 1 in [6] with the same conditions.

It is a small nuisance that (12) is not as explicit as condition (21) in the i.i.d. case, as (12) is given in terms of the values of $J_1$ and $J_2$, which in turn are the results of an optimization. We showed above how to solve this optimization problem explicitly in the i.i.d. case, but it does not seem that there exists such a simple solution for general Markov chains. From a practical point of view one may notice that taking $g^*(x, y, x', y') = 3\theta^* f(x', y')/4$, then

(22) $$\max\{\varphi_1(g^*), \varphi_2(g^*)\} < 1$$

implies (12). Since $\varphi_1(g^*)$ and $\varphi_2(g^*)$ can be computed numerically we see that (22) provides a usable, sufficient criterion for Theorem 3.1 to hold.

**4. The counting construction.** We will show that $C_n(t_n)$ is asymptotically Poisson distributed by constructing another counting variable, which equals $C_n(t_n)$ with probability tending to 1, and for which we can verify the conditions given in Theorem 1 in [2].

We need to introduce some notation. Let

$$I = \{(i,j) | 0 \leq i, j \leq n-1\};$$

then for each $a = (i,j) \in I$ and $\delta > 0$ we define the (pair) empirical measure $\varepsilon_{a,\delta}$ by

$$\varepsilon_{a,\delta}((x,y),(x',y')) = \frac{1}{\delta}\sum_{k=1}^{\delta} \mathbf{1}_{(x,y),(x',y')}((X_{i+k-1}, Y_{j+k-1}),(X_{i+k}, Y_{j+k}))$$

for $(x,y), (x', y') \in E^2$. With abuse of notation we will in the following also use $f$ to denote the function defined on $E^2 \times E^2$ by $(x, y, x', y') \mapsto f(x', y')$. Then

$$\delta\varepsilon_{a,\delta}(f) = \sum_{k=1}^{\delta} f(X_{i+k}, Y_{j+k}) = S_{i,j}^\delta.$$

Let $d$ denote the total variation metric on the set of probability measures on $E^2 \times E^2$. Then for $a \in I$ and for any $t > 0$, $\eta > 0$ and integer $l > 0$ define the variable

$$V_a = V_a(t, l, \eta) = 1\left(T_a = 0, \max_{\substack{\delta:\, \delta \leq \Delta(a) \wedge l \\ \text{and } d(\varepsilon_{a,\delta}, \hat{\pi}) < \eta}} \delta\varepsilon_{a,\delta}(f) > t\right).$$

We should observe that the counting variable $C_n(t)$ has the following representation:

(23) $$C_n(t) = \sum_{a \in I} 1\left(T_a = 0, \max_{\delta:\, \delta \leq \Delta(a)} \delta\varepsilon_{a,\delta}(f) > t\right).$$

We show in Section 5.8 that in the setup of the present paper, for a suitable choice of $l_n$ and $\eta > 0$, then

(24) $$\mathbb{P}\left(\sum_{a \in I} V_a(t_n, l_n, \eta) \neq C_n(t_n)\right) \to 0$$

when $n \to \infty$. The reason for introducing the $l$-restriction is to be able to control the dependencies between the $V_a$-variables better. The reason for the restriction on the empirical measures is more subtle, and we give a discussion of this in Section 6.

As mentioned, we prove that $\sum_{a \in I} V_a$ is asymptotically Poisson distributed by applying Theorem 1 in [2], which is based on the Chen–Stein method.

We assume that a subset $B_a \subseteq I$ is given for all $a \in I$. This set $B_a$ is called the neighborhood of strong dependence of $V_a$, and in the proof of Lemma 5.16 we make a concrete choice of $B_a$. Furthermore, for $a \in I$ let

$$\mathcal{F}_a = \sigma(V_b | b \notin B_a)$$

denote the $\sigma$-algebra generated by those variables $V_b$ *not* in the neighborhood of strong dependence of $V_a$.

Rephrasing Theorem 1 in [2] gives:

THEOREM 4.1. *Suppose that $(l_n)_{n \geq 1}$, $(t_n)_{n \geq 1}$ and $\eta > 0$ are chosen such that for some sequence $(\lambda_n)_{n \geq 1}$*

(25) $$\beta_{1,n} = \left|\sum_{a \in I} \mathbb{E}(V_a) - \lambda_n\right| \to 0,$$

*for $n \to \infty$, and suppose that*

(26) $$\beta_{2,n} = \sum_{a \in I, b \in B_a} \mathbb{E}(V_a)\mathbb{E}(V_b) \to 0,$$

(27) $$\beta_{3,n} = \sum_{a \in I, b \in B_a, b \neq a} \mathbb{E}(V_a V_b) \to 0,$$

(28) $$\beta_{4,n} = \sum_{a \in I} \mathbb{E}|\mathbb{E}(V_a | \mathcal{F}_a) - \mathbb{E}(V_a)| \to 0,$$

*for* $n \to \infty$; *then*

$$\left\| \mathcal{D}\left(\sum_{a \in I} V_a\right) - \mathrm{Poi}(\lambda_n) \right\| \to 0. \tag{29}$$

*In fact, the bound*

$$\left\| \mathcal{D}\left(\sum_{a \in I} V_a\right) - \mathrm{Poi}(\lambda_n) \right\| \leq \beta_{1,n} + 2(\beta_{2,n} + \beta_{3,n} + \beta_{4,n})$$

*always holds.*

As a direct consequence, using the coupling inequality, we have the following corollary.

COROLLARY 4.2. *If* (29) *holds and* (24) *is fulfilled also, then*

$$\|\mathcal{D}(C_n(t_n)) - \mathrm{Poi}(\lambda_n)\| \to 0 \tag{30}$$

*and*

$$\mathbb{P}(\mathcal{M}_n \leq t_n) - \exp(-\lambda_n) \to 0. \tag{31}$$

**5. Proofs.** The proof of Theorem 3.1 is divided into a number of lemmas. We need to verify the conditions in Theorem 4.1, and to this end we need bounds on the expectations $\mathbb{E}(V_a V_b) = \mathbb{P}(V_a = 1, V_b = 1)$ for $b \in B_a$ and $a \neq b$. This is the subject of the following subsections and the most difficult part of the proof. In Section 5.8 we collect the bounds obtained to prove that the conditions of Theorem 4.1 are fulfilled when aligning independent Markov chains under the assumptions given in Theorem 3.1 and we show that (24) holds.

For $a, b \in I$ we always have that

$$\begin{aligned}
\mathbb{E}(V_a V_b) &\leq \mathbb{P}\left( \max_{\substack{\delta : \delta \leq \Delta(a) \wedge l \\ \text{and } d(\varepsilon_{a,\delta}, \hat{\pi}) < \eta}} \delta \varepsilon_{a,\delta}(f) > t, \max_{\substack{\delta : \delta \leq \Delta(b) \wedge l \\ \text{and } d(\varepsilon_{b,\delta}, \hat{\pi}) < \eta}} \delta \varepsilon_{b,\delta}(f) > t \right) \\
&\leq l^2 \max_{1 \leq \delta_1, \delta_2 \leq l} \mathbb{P}\left( \begin{array}{l} \delta_1 \varepsilon_{a,\delta_1}(f) > t, d(\varepsilon_{a,\delta_1}, \hat{\pi}) < \eta, \\ \delta_2 \varepsilon_{b,\delta_2}(f) > t, d(\varepsilon_{b,\delta_2}, \hat{\pi}) < \eta \end{array} \right).
\end{aligned} \tag{32}$$

To bound $\mathbb{E}(V_a V_b)$ we thus need to bound the probability on the right-hand side above. The same $X$- and $Y$-variables may enter both of the empirical measures in two essentially different ways. Either variables from both sequences enter both empirical measures or only variables from one sequence enter both empirical measures. These two different cases need different treatment. To give an exhaustive treatment of the different ways that such a sharing of variables can be arranged becomes unreasonably complicated, so we choose to treat the two essentially different cases for a specific arrangement of the sharing of variables in sufficient detail for the reader to be able to convince himself that all other arrangements can be treated similarly.

5.1. *Positive functionals of a Markov chain.* We make a useful and general observation on how to bound the expectation of positive functionals of a Markov chain. It allows us to assume parts of the same Markov chain to be independent, stationary versions at the expense of a constant factor.

LEMMA 5.1. *Let $\mathcal{Z} = (Z_n)_{n \geq 0}$ be an irreducible Markov chain on a finite state space $F$ and let $0 = k_1 < \cdots < k_N < \infty$ be given. Then there exists a constant $\rho_N$ such that if $(Z_n^i)_{n=k_i}^{k_{i+1}}$ for $i = 1, \ldots, N$ ($k_{N+1} = \infty$) are $N$ independent stationary Markov chains with the same transition probabilities as $\mathcal{Z}$, and $\widetilde{\mathcal{Z}} = (\widetilde{Z}_n)_{n \geq 0}$ is given by $\widetilde{Z}_n = Z_n^i$ if $k_i \leq n < k_{i+1}$, then for a positive functional*

$$\Lambda : F^{\mathbb{N}_0} \to [0, \infty)$$

*it holds that*

$$\mathbb{E}(\Lambda(\mathcal{Z})) \leq \rho_N \mathbb{E}(\Lambda(\widetilde{\mathcal{Z}})). \tag{33}$$

*The constant $\rho_N$ does not depend on the actual initial distribution of $\mathcal{Z}$ nor on the functional $\Lambda$.*

PROOF. Assume $N = 2$. The general result follows by induction. Assume first that $\mathcal{Z}$ is stationary and that $(Z_n^1)_{n=0}^{k_2}$ and $(Z_n^2)_{n \geq k_2}$ are independent and stationary. Then $\mathcal{Z}$ has the same distribution as $\widetilde{\mathcal{Z}}$ conditionally on $Z_{k_2}^1 = Z_{k_2}^2$; hence using that $\Lambda$ is a positive functional

$$\mathbb{E}(\Lambda(\mathcal{Z})) = \frac{\mathbb{E}(\Lambda(\widetilde{\mathcal{Z}}); Z_{k_2}^1 = Z_{k_2}^2)}{\mathbb{P}(Z_{k_2}^1 = Z_{k_2}^2)}$$

$$\leq \rho \mathbb{E}(\Lambda(\widetilde{\mathcal{Z}}))$$

with $\rho = (\sum_{x \in E} \pi_x^2)^{-1}$, where $\pi$ is the invariant distribution.

If $\mathcal{Z}$ is nonstationary with initial distribution $\nu$, say, we have that

$$\mathbb{E}_\nu(\Lambda(\mathcal{Z})) = \sum_x \frac{\nu_x}{\pi_x} \pi_x \mathbb{E}_x(\Lambda(\mathcal{Z}))$$

$$\leq \frac{1}{\min_x \pi_x} \mathbb{E}_\pi(\Lambda(\mathcal{Z})).$$

So $\rho_2 = \rho / \min_x \pi_x$ will do. In general $\rho_N = \rho^{N-1} / \min_x \pi_x$ can be used. □

5.2. *Exponential change of measure.* Let $\mathcal{Z} = (Z_n)_{n \geq 0}$ be a Markov chain on a finite state space $F$ with transition probabilities $R$. Assume that $R$ is irreducible, and assume that $g : F \times F \to \mathbb{R}$ is a given function. Then we define the matrix $\Psi(g)$ with positive entries by

$$\Psi(g)_{x,x'} = \exp(g(x,x'))R_{x,x'},$$

with spectral radius $\psi(g)$ and corresponding right eigenvector $r^g = (r^g(x))_{x \in F}$. Due to irreducibility of $\Psi(g)$ this eigenvector has strictly positive entries. With

$$g_n(\mathcal{Z}) = \sum_{k=1}^{n} g(Z_{k-1}, Z_k)$$

we define the process $(L_n^g)_{n \geq 0}$ by

$$L_n^g = \frac{r^g(Z_n) \exp(g_n(\mathcal{Z}))}{r^g(Z_0) \psi(g)^n}.$$

Then with $(\mathcal{F}_n)_{n \geq 0}$ the filtration of $\sigma$-algebras generated by the Markov chain it follows that

$$\mathbb{E}(L_n^g | \mathcal{F}_{n-1}) = \frac{\exp(g_{n-1}(\mathcal{Z}))}{r^g(Z_0)\psi(g)^n} \mathbb{E}(r^g(Z_n) \exp(g(Z_{n-1}, Z_n)) | Z_{n-1})$$

(34)
$$= \frac{\exp(g_{n-1}(\mathcal{Z}))}{r^g(Z_0)\psi(g)^n} (\Psi(g)r^g)(Z_{n-1})$$

$$= \frac{\exp(g_{n-1}(\mathcal{Z}))}{r^g(Z_0)\psi(g)^n} \psi(g)r^g(Z_{n-1}) = L_{n-1}^g.$$

This shows that $(L_n^g, \mathcal{F}_n)_{n \geq 0}$ is a martingale, for which $L_n^g > 0$ and $\mathbb{E}(L_n^g) = \mathbb{E}(L_0^g) = 1$. A probability measure $\mathbb{P}_n^g$ on $\mathcal{F}_n$ is then defined to have Radon–Nikodym derivative $L_n^g$ w.r.t. the restriction of $\mathbb{P}$ to $\mathcal{F}_n$. These measures can be extended to a single measure $\mathbb{P}^g$, the exponentially changed or exponentially tilted measure, under which $(Z_n)_{n \geq 0}$ is a Markov chain ([3], Theorem XIII.8.1) with transition probabilities

$$R_{x,x'}^g = \frac{r^g(x')}{r^g(x)\psi(g)} \Psi(g)_{x,x'}.$$

We should observe that since the eigenvector fraction is bounded below by a strictly positive constant, and bounded above as well, then $\mathbb{E}(L_n^g) = 1$ implies that

(35) $$\frac{1}{n} \log \mathbb{E}(\exp(g_n(\mathcal{Z}))) \to \log \psi(g)$$

for $n \to \infty$.

If we return to the setup of the present paper with $F = E^2$, $g = \theta f$ for $\theta \in \mathbb{R}$, and the Markov chain being $\mathcal{Z} = (X_n, Y_n)_{n \geq 0}$, then

$$g_n(\mathcal{Z}) = \sum_{k=1}^{n} \theta f(X_k, Y_k) = \theta S_n,$$

and we find that the matrix $R^*$ introduced in Section 3 is precisely the matrix of transition probabilities for the Markov chain under the exponentially changed measure $\mathbb{P}^{\theta^* f}$. We will denote this measure simply by $\mathbb{P}^*$. Note that the exponential change of measure does not change the distribution $\pi$ of $(X_0, Y_0)$ whereas the invariant measure $\pi^*$ for $R^*$ typically differs from $\pi$. The measure under which $(X_n, Y_n)_{n \geq 0}$ is a *stationary* Markov chain with transition probabilities $R^*$ will be denoted $\mathbb{P}^*_{\pi^*}$. We use $\mathbb{E}^*$ to denote expectations under $\mathbb{P}^*$.

If we define the stopping time $\tau = \inf\{n \geq 0 | S_n > t\}$, then an easy consequence of the exponential change of measure technique is, according to [3], Theorem XIII.3.2, the following *Lundberg-type inequality*: For any event $G \in \mathcal{F}_\tau$ with $G \subseteq (\tau < \infty)$

$$(36) \qquad \mathbb{P}(G) = \mathbb{E}^*\left(\frac{1}{L^*_\tau}; G\right) \leq K \exp(-\theta^* t).$$

The inequality follows from $L^g_\tau \geq K \exp(\theta^* t)$ where $K$ bounds that eigenvector fraction.

5.3. *Variables shared in one sequence.* Let $g : E^2 \times E^2 \to \mathbb{R}$ be a function and let $r_i^g = (r_i^g(x, y, z))$ denote the right eigenvector for $\Phi_i(g)$ with eigenvalue $\varphi_i(g)$ for $i = 1, 2$, respectively. As above, due to irreducibility, all coordinates of these vectors are strictly positive.

In this section we derive a result corresponding to variables shared from the $X$-sequence only, and we thus use the $\Phi_1$ matrix. Similar derivations for variables shared from the $Y$-sequence only using $\Phi_2$ are possible.

Fix $i \leq \delta_1$ and $\delta_2 \geq \delta_1 - i$ and define the functions

$$\sigma_1((x_k)_k, (y_k)_k) = \sum_{k=1}^{i} g(x_{k-1}, y_{k-1}, x_k, y_k),$$

$$\sigma_2((x_k)_k, (y_k)_k, (z_k)_k) = \sum_{k=i+1}^{\delta_1} g(x_{k-1}, y_{k-1}, x_k, y_k) + g(x_{k-1}, z_{k-1}, x_k, z_k),$$

$$\sigma_3((x_k)_k, (z_k)_k) = \sum_{k=\delta_1+1}^{i+\delta_2} g(x_{k-1}, z_{k-1}, x_k, z_k).$$

Let $\Phi_0(g)$ denote the matrix

$$\Phi_0(g)_{(x,y),(x',y')} = \exp(g(x, y, x', y')) P_{x,x'} Q_{y,y'}$$

and $\varphi_0(g)$ the spectral radius. Let $r_0^g$ denote the corresponding right eigenvector. We define a positive functional $\mathcal{L}^g$ on $(E^3)^{\mathbb{N}_0}$ by

$$\mathcal{L}^g = \frac{r_0^g(x_i, y_i) \exp(\sigma_1)}{r_0^g(x_0, y_0) \varphi_0(g)^i} \frac{r_1^g(x_{\delta_1}, y_{\delta_1}, z_{\delta_1}) \exp(\sigma_2)}{r_1^g(x_i, y_i, z_i) \varphi_1(g)^{\delta_1 - i}} \frac{r_0^g(x_{i+\delta_2}, z_{i+\delta_2}) \exp(\sigma_3)}{r_0^g(x_{\delta_1}, z_{\delta_1}) \varphi_0(g)^{i+\delta_2-\delta_1}}.$$

Assume that $(Z_n)_{n \geq 1}$ is a stationary Markov chain with transition probabilities $Q$ independent of $(X_n, Y_n)_{n \geq 1}$, and let $\mathcal{X} = (X_n)_{n \geq 1}$, $\mathcal{Y} = (Y_n)_{n \geq 1}$ and $\mathcal{Z} = (Z_n)_{n \geq 1}$. Introduce also $\mathcal{Y}^T = (Y_{T+n})_{n \geq 1}$ as the $T$-shift of $\mathcal{Y}$ for $T \geq 1$.

LEMMA 5.2. *It holds that*

$$\mathbb{E}(\mathcal{L}^g(\mathcal{X}, \mathcal{Y}, \mathcal{Z})) = 1, \tag{37}$$

*and, furthermore, there exists a constant $\rho > 0$ such that*

$$\mathbb{E}(\mathcal{L}^g(\mathcal{X}, \mathcal{Y}, \mathcal{Y}^T)) \leq \rho \tag{38}$$

*whenever $i + T \geq \delta_1 + 1$.*

PROOF. The first part of the lemma follows by repeating the arguments in (34) three times corresponding to making three different, successive exponential changes of measures. The second claim follows by Lemma 5.1. □

We restrict our attention to the case where $i + T \geq \delta_1 + 1$, so that there is no overlap in the $Y$-sequence. Let $\varepsilon_1 = \varepsilon_{(0,0),\delta_1}$ and $\varepsilon_2 = \varepsilon_{(i,i+T),\delta_2}$.

LEMMA 5.3. *For any $g : E \times E \to \mathbb{R}$ and $\varepsilon > 0$ there exist constants $\eta, K > 0$ such that for all $s > 0$,*

$$\mathbb{P}\begin{pmatrix} \delta_1 \varepsilon_1(f) > s, d(\varepsilon_1, \hat{\pi}) < \eta, \\ \delta_2 \varepsilon_2(f) > s, d(\varepsilon_2, \hat{\pi}) < \eta \end{pmatrix} \leq K \exp\left(-s\left(\frac{2\hat{\pi}(g) - \log \varphi_1(g)}{\pi^*(f)} - \varepsilon\right)\right).$$

PROOF. First we show that $\log \varphi_1(g) \geq 2 \log \varphi_0(g)$. Let

$$g_n(\mathcal{X}, \mathcal{Y}) = \sum_{k=1}^{n} g(X_{k-1}, Y_{k-1}, X_k, Y_k)$$

and define $g_n(\mathcal{X}, \mathcal{Z})$ likewise. Let $\mathbb{E}_\mathcal{X}$, $\mathbb{E}_\mathcal{Y}$ and $\mathbb{E}_\mathcal{Z}$ denote the expectation operators where we only integrate w.r.t. the distribution of $\mathcal{X}$, $\mathcal{Y}$ or $\mathcal{Z}$, respectively. Introduce

$$\rho_n(\mathcal{X}) = \mathbb{E}_\mathcal{Y}(\exp(g_n(\mathcal{X}, \mathcal{Y})));$$

then by (35) and Tonelli

$$\log \varphi_0(g) = \lim_{n \to \infty} \frac{1}{n} \log \mathbb{E}_{\mathcal{X}}(\rho_n(\mathcal{X})).$$

Using Tonelli again and Jensen's inequality, and that $\mathcal{Y}$ and $\mathcal{Z}$ are independent and identically distributed, we find that

$$\mathbb{E}(\exp(g_n(\mathcal{X}, \mathcal{Y}) + g_n(\mathcal{X}, \mathcal{Z}))) = \mathbb{E}_{\mathcal{X}}(\mathbb{E}_{\mathcal{Y}}(\exp(g_n(\mathcal{X}, \mathcal{Y})))\mathbb{E}_{\mathcal{Z}}(\exp(g_n(\mathcal{X}, \mathcal{Z}))))$$
$$= \mathbb{E}_{\mathcal{X}}(\rho_n(\mathcal{X})^2) \geq \mathbb{E}_{\mathcal{X}}(\rho_n(\mathcal{X}))^2.$$

Using (35) again gives

$$\log \varphi_1(g) = \lim_{n \to \infty} \frac{1}{n} \log \mathbb{E}(\exp(g_n(\mathcal{X}, \mathcal{Y}) + g_n(\mathcal{X}, \mathcal{Z})))$$
$$\geq 2 \lim_{n \to \infty} \frac{1}{n} \log \mathbb{E}_{\mathcal{X}}(\rho_n(\mathcal{X})) = 2 \log \varphi_0(g).$$

Since $2(\delta_1 - i) + i + (i + \delta_2 - \delta_1) = \delta_1 + \delta_2$ and $\sigma_1 + \sigma_2 + \sigma_3 = \delta_1 \varepsilon_1(g) + \delta_2 \varepsilon_2(g)$, the inequality $\log \varphi_1(g) \geq 2 \log \varphi_0(g)$ gives that

$$\mathcal{L}^g(\mathcal{X}, \mathcal{Y}, \mathcal{Y}^T) \geq \gamma \exp(\delta_1 \varepsilon_1(g) + \delta_2 \varepsilon_2(g) - (\delta_1 + \delta_2) \log \varphi_1(g)/2)$$

with $\gamma > 0$ a lower bound on the eigenvector fractions. We may assume that $2\hat{\pi}(g) - \log \varphi_1(g) > \pi^*(f)\varepsilon$ since the result is trivial otherwise. Then we can find $\varepsilon' > 0$ such that

$$\frac{2(\hat{\pi}(g) - \varepsilon') - \log \varphi_1(g)}{\pi^*(f) + \varepsilon'} = \frac{2\hat{\pi}(g) - \log \varphi_1(g)}{\pi^*(f)} - \varepsilon$$

and choose $\eta$ so small that for $\nu$ a probability measure on $E^2 \times E^2$ with $d(\nu, \hat{\pi}) < \eta$ we have $|\nu(g) - \hat{\pi}(g)| \leq \varepsilon'$ and $|\nu(f) - \pi^*(f)| \leq \varepsilon'$. On the event

$$A = \begin{pmatrix} \delta_1 \varepsilon_1(f) > s, d(\varepsilon_1, \hat{\pi}) < \eta, \\ \delta_2 \varepsilon_2(f) > s, d(\varepsilon_2, \hat{\pi}) < \eta \end{pmatrix}$$

we see that

$$\delta_1 \varepsilon_1(g) + \delta_2 \varepsilon_2(g) - (\delta_1 + \delta_2) \log \varphi_1(g)/2 \geq \frac{\delta_1 + \delta_2}{2}(2(\hat{\pi}(g) - \varepsilon') - \log \varphi_1(g))$$
$$\geq s\left(\frac{2\hat{\pi}(g) - \log \varphi_1(g)}{\pi^*(f)} - \varepsilon\right)$$

since on $A$ we have $\delta_1 + \delta_2 > 2s/(\pi^*(f) + \varepsilon')$. Hence

$$\mathbb{P}(A) = \mathbb{E}\left(\frac{\mathcal{L}^g(\mathcal{X}, \mathcal{Y}, \mathcal{Y}^T)}{\mathcal{L}^g(\mathcal{X}, \mathcal{Y}, \mathcal{Y}^T)}; A\right)$$
$$\leq \gamma^{-1} \exp\left(-s\left(\frac{2\hat{\pi}(g) - \log \varphi_1(g)}{\pi^*(f)} - \varepsilon\right)\right)\mathbb{E}(\mathcal{L}^g(\mathcal{X}, \mathcal{Y}, \mathcal{Y}^T); A)$$
$$\leq \rho \gamma^{-1} \exp\left(-s\left(\frac{2\hat{\pi}(g) - \log \varphi_1(g)}{\pi^*(f)} - \varepsilon\right)\right),$$

where the first inequality follows by bounding the denominator from below using the inequalities obtained above, and the second inequality follows from Lemma 5.2. □

If the condition (12) is fulfilled, then $2J_1 > 3\theta^*\pi^*(f)$ and we can in particular choose a function $g$ and an $\varepsilon > 0$ sufficiently small such that

$$2\hat{\pi}(g) - \log \varphi_1(g) \geq (3\theta^*/2 + 2\varepsilon)\pi^*(f).$$

The following corollary is therefore a direct consequence of Lemma 5.3.

COROLLARY 5.4. *If* (12) *is fulfilled, there exist constants* $\varepsilon, \eta, K > 0$ *such that for all* $s > 0$

$$\tag{39} \mathbb{P}\begin{pmatrix} \delta_1\varepsilon_1(f) > s, d(\varepsilon_1, \hat{\pi}) < \eta, \\ \delta_2\varepsilon_2(f) > s, d(\varepsilon_2, \hat{\pi}) < \eta \end{pmatrix} \leq K\exp(-(3\theta^*/2 + \varepsilon)s).$$

The result in Corollary 5.4 gives a prototypical inequality under the assumption $2J_1 > 3\theta^*\pi^*(f)$ when only variables from the $X$-sequence enter both of the empirical measures. If only variables from the $Y$-sequence enter both empirical measures, a similar inequality is obtained under the assumption that $2J_2 > 3\theta^*\pi^*(f)$.

5.4. *A uniform large deviation result.* To handle the case with variables shared from both sequences we need a special large deviation result for Markov chains that we will derive in this section. We first state the useful Azuma–Hoeffding inequality for martingales with bounded increments; see Lemma 1.5 in [13] or Theorem 1.3.1 in [19].

LEMMA 5.5. *If* $(Z_n, \mathcal{F}_n)_{n \geq 0}$ *is a mean-zero martingale with* $Z_0 = 0$ *such that for all* $n \geq 1$

$$|Z_n - Z_{n-1}| \leq c_n$$

*for some sequence* $(c_n)_{n \geq 1}$, *then for* $\lambda > 0$

$$\mathbb{P}(Z_n \geq \lambda) \leq \exp\left(-\frac{\lambda^2}{2\sum_{k=1}^n c_k^2}\right).$$

Fix $j \geq 1$ and let in this section $(X_n, Y_n)_{n=1}^j$ be a stationary, aperiodic and irreducible Markov chain with transition probabilities given by $R$ and invariant distribution $\pi_R$. Let $(Y_n)_{n \geq j+1}$ be an independent, stationary, aperiodic and irreducible Markov chain with transition probabilities given by $Q$ and invariant distribution $\pi_Q$. For an $E^2 \times E^2$ matrix $G$ define the norm of the matrix as

$$\|G\|_\infty = \max_{(x,y)} \sum_{(z,w)} |G_{(x,y),(z,w)}|.$$

With $\mathbb{1}$ the column vector of 1's, the matrix $R^k$ converges to $\mathbb{1}\pi_R$ due to irreducibility and aperiodicity, and since the rate of convergence is sufficiently fast, in fact geometric, we have that

$$\sum_{k=0}^{\infty} \|R^k - \mathbb{1}\pi_R\|_\infty < \infty.$$

For an $E^2$ vector $f$ we let $\|f\|_\infty = \max_{(x,y)} |f(x,y)|$ denote the max-norm. Then clearly for any $E^2 \times E^2$ matrix $G$, with $G(f)$ the matrix product of $G$ with the vector $f$, $\|G(f)\|_\infty \leq \|f\|_\infty \|G\|_\infty$, and especially

$$\|R^k(f) - \mathbb{1}\pi_R(f)\|_\infty \leq \|f\|_\infty \|R^k - \mathbb{1}\pi_R\|_\infty.$$

For $T \geq 1$ a fixed constant we want to give an exponential bound of the probability

(40) $$\mathbb{P}\left(\sum_{k=1}^{j} f(X_k, Y_{k+T}) \geq \sum_{k=1}^{j} f(X_k, Y_k)\right)$$

if $\mathbb{E}(f(X_k, Y_{k+T})) < \mathbb{E}(f(X_k, Y_k))$ all $k$. This is achieved by introducing a relevant martingale and then using the Azuma–Hoeffding inequality.

Let $\mathcal{F}_0 = \{\varnothing, \Omega\}$ and for $n \geq 1$ let $\mathcal{F}_n$ denote the $\sigma$-algebra generated by $X_1, \ldots, X_n, Y_1, \ldots, Y_n$ together with $Y_{j+1}, \ldots, Y_{n+T}$ if $n + T > j$. Define

$$S_{j,T} = \sum_{k=1}^{j} [f(X_k, Y_{k+T}) - f(X_k, Y_k)] \qquad (S_{0,T} = 0),$$

and with $\xi_{j,T} = \mathbb{E}(S_{j,T})$ let

(41) $$Z_n = \mathbb{E}(S_{j,T} - \xi_{j,T} | \mathcal{F}_n).$$

Then $(Z_n, \mathcal{F}_n)_{n=0}^{j}$ is a mean-zero martingale with $Z_0 = 0$ (depending on $T$, though we have suppressed this in the notation). Notice that $Z_j = S_{j,T} - \xi_{j,T}$. The following lemma shows that the martingale differences

$$|Z_n - Z_{n-1}| = |\mathbb{E}(S_{j,T} | \mathcal{F}_n) - \mathbb{E}(S_{j,T} | \mathcal{F}_{n-1})|$$

are uniformly bounded by a constant.

LEMMA 5.6. *There exists a constant $\eta$ independent of $j$ and $T$ such that*

(42) $$|Z_n - Z_{n-1}| \leq \eta.$$

*Here $\eta$ can be chosen as*

$$\eta = 6\|f\|_\infty \sum_{k=0}^{\infty} \|R^k - \mathbb{1}\pi_R\|_\infty. \tag{43}$$

PROOF. The Markov property gives that for $n \leq k \leq j$

$$\mathbb{E}(f(X_k, Y_k)|\mathcal{F}_n) = R^{k-n}(f)(X_n, Y_n).$$

Define the function $\hat{f}$ by

$$\hat{f}(x,y) = R^T(f(x,\cdot))(x,y) = \sum_{z,w} f(x,w) R^T_{(x,y),(z,w)},$$

and for $n \leq k \leq j$ define $\tilde{f}_{k,n}$

$$\tilde{f}_{k,n}(x,y) = \begin{cases} \sum_z f(x,z) Q^{k-n}_{y,z}, & \text{if } n+T > j, \\ \sum_z f(x,z) \pi_Q(z), & \text{if } n+T \leq j. \end{cases}$$

Then

$$\mathbb{E}(f(X_k, Y_{k+T})|\mathcal{F}_n) = \begin{cases} R^{k-n}(\hat{f})(X_n, Y_n), & k \in C_1, \\ R^{k-n}(\tilde{f}_{k,n}(\cdot, Y_{n+T}))(X_n, Y_n), & k \in C_2, \\ R^{k+T-n}(f(X_k, \cdot))(X_n, Y_n), & k \in C_3, \end{cases}$$

where

$$C_1 = \{k | n \leq k < k+T \leq j\},$$
$$C_2 = \{k | n \leq k \leq j < k+T\},$$
$$C_3 = \{k | n-T \leq k < n \leq k+T \leq j\}.$$

Observing that

$$\mathbb{E}(S_{j,T}|\mathcal{F}_n) = \sum_{k=1}^{j} \mathbb{E}(f(X_k, Y_{k+T})|\mathcal{F}_n) - \sum_{k=1}^{j} \mathbb{E}(f(X_k, Y_k)|\mathcal{F}_n)$$

and subtracting $\mathbb{E}(S_{j,T}|\mathcal{F}_{n-1})$ from this, the martingale difference $Z_n - Z_{n-1}$ is seen to be the sum of the following two terms:

$$t_1 = \sum_{k=n-T}^{j} [\mathbb{E}(f(X_k, Y_{k+T})|\mathcal{F}_n) - \mathbb{E}(f(X_k, Y_{k+T})|\mathcal{F}_{n-1})],$$

$$t_2 = \sum_{k=n}^{j} [\mathbb{E}(f(X_k, Y_k)|\mathcal{F}_{n-1}) - \mathbb{E}(f(X_k, Y_k)|\mathcal{F}_n)].$$

Since
$$|\mathbb{E}(f(X_k,Y_k)|\mathcal{F}_n) - \pi_R(f)| = |R^{k-n}(f)(X_n,Y_n) - \pi_R(f)|$$
$$\leq \|f\|_\infty \|R^{k-n} - \mathbb{1}\pi_R\|_\infty,$$

the term $t_2$ is controlled by the following inequality:

$$(44) \quad |t_2| \leq 2\|f\|_\infty \sum_{k=n}^{j} \|R^{k-n} - \mathbb{1}\pi_R\|_\infty \leq 2\|f\|_\infty \sum_{k=0}^{\infty} \|R^k - \mathbb{1}\pi_R\|_\infty.$$

Noting that $\|\hat{f}\|_\infty, \|\tilde{f}_{k,n}(\cdot,y)\|_\infty, \|f(x,\cdot)\|_\infty \leq \|f\|_\infty$ we observe that for $k \in C_1$,

$$|\mathbb{E}(f(X_k,Y_{k+T})|\mathcal{F}_n) - \pi_R(\hat{f})| \leq \|f\|_\infty \|R^{k-n} - \mathbb{1}\pi_R\|_\infty,$$

for $k \in C_2$,

$$|\mathbb{E}(f(X_k,Y_{k+T})|\mathcal{F}_n) - \pi_R(\tilde{f}_{k,n}(\cdot,Y_{n+T}))| \leq \|f\|_\infty \|R^{k-n} - \mathbb{1}\pi_R\|_\infty,$$

and for $k \in C_3$,

$$|\mathbb{E}(f(X_k,Y_{k+T})|\mathcal{F}_n) - \pi_R(f(X_k,\cdot))| \leq \|f\|_\infty \|R^{k+T-n} - \mathbb{1}\pi_R\|_\infty.$$

Since the three inequalities above also hold when conditioning on $\mathcal{F}_{n-1}$ we obtain

$$\sum_{k \in C_1 \cup C_2 \cup C_3} |\mathbb{E}(f(X_k,Y_{k+T})|\mathcal{F}_n) - \mathbb{E}(f(X_k,Y_{k+T})|\mathcal{F}_{n-1})|$$
$$\leq 2\|f\|_\infty \sum_{k \in C_1 \cup C_2} \|R^{k-n} - \mathbb{1}\pi_R\|_\infty + 2\|f\|_\infty \sum_{k \in C_3} \|R^{k+T-n} - \mathbb{1}\pi_R\|_\infty$$
$$\leq 4\|f\|_\infty \sum_{k=0}^{\infty} \|R^k - \mathbb{1}\pi_R\|_\infty.$$

Finally, if $n - T \leq k < n < j < k + T$, then

$$\mathbb{E}(f(X_k,Y_{k+T})|\mathcal{F}_n) = \mathbb{E}(f(X_k,Y_{k+T})|\mathcal{F}_{n-1}) = f(X_k,Y_{k+T}),$$

hence

$$|t_1| \leq 4\|f\|_\infty \sum_{k=0}^{\infty} \|R^k - \mathbb{1}\pi_R\|_\infty,$$

which together with (44) gives (42) with $\eta$ chosen as (43). □

THEOREM 5.7. *If $\xi_{j,T} < 0$, it holds that*

$$(45) \quad \mathbb{P}(S_{j,T} \geq 0) = \mathbb{P}(S_{j,T} - \xi_{j,T} \geq -\xi_{j,T}) \leq \exp\left(-\frac{\xi_{j,T}^2}{2j\eta^2}\right)$$

*with $\eta$ chosen as in Lemma 5.6.*

PROOF. This follows directly from the Azuma–Hoeffding inequality for the mean-zero martingale $(Z_n, \mathcal{F}_n)_{n=1}^{j}$, since it has increments uniformly bounded by $\eta$. □

5.5. *Mean value inequalities.* We will apply the result in the previous section by considering the Markov chain $(X_n, Y_n)_{n=1}^{j}$ under the exponentially tilted measure $\mathbb{P}_{\pi^*}^*$ and $(Y_n)_{n \geq j+1}$ under $\mathbb{P}_\pi$. To do so we will need to establish inequalities relating the mean of $f(X_n, Y_n)$ to the mean of $f(X_n, Y_{n+T})$ [or $f(X_{n+T}, Y_n)$]. Let

$$\mu^* = \mathbb{E}_{\pi^*}^*(f(X_1, Y_1)) = \pi^*(f)$$

denote the stationary mean of $f(X_n, Y_n)$ under the exponentially tilted measure and let

$$\mu_T^* = \mathbb{E}_{\pi^*}^*(f(X_1, Y_{1+T}))$$

denote the stationary mean when shifting the $Y$-sequence $T$ positions.

It was mentioned in Section 3 that the function $\varphi$ is log-convex. In the following we will need results obtained in [16] about *strict* log-convexity of $\varphi$-like functions.

Let $F$ be a finite set, $g: F \to \mathbb{R}$ any function and $R$ an irreducible $F \times F$ matrix of transition probabilities. Following Definition 2 in [16] we say that $g$ is degenerate w.r.t. $R$ if there exists a constant $\gamma$ such that for all cycles $x_1, \ldots, x_n$ w.r.t. $R$

$$\sum_{k=1}^{n} g(x_k) = \gamma n.$$

Let $\psi(\theta)$ for $\theta \in \mathbb{R}$ be the spectral radius of the $F \times F$ matrix $\Psi(\theta)$ with entries

$$\Psi(\theta)_{x,x'} = \exp(\theta g(x')) R_{x,x'}.$$

From Theorem 5 in [16] it follows that if $g$ is nondegenerate w.r.t. $R$, then $\psi$ is strictly log-convex, and if $g$ is degenerate w.r.t. $R$, then $\psi(\theta) = \exp(\gamma \theta)$ (i.e., $\log \psi$ is linear). The consequence that we will use repeatedly below is that if $\psi(0) = \psi(\theta^*) = 1$ for $\theta^* > 0$ and if $g$ is degenerate w.r.t. $R$, then necessarily $\psi(\theta) = 1$ for all $\theta \in \mathbb{R}$ and the constant $\gamma$ equals 0. Thus if we can find a single cycle $x_1, \ldots, x_n$ w.r.t. $R$ such that

$$\sum_{k=1}^{n} g(x_k) \neq 0, \tag{46}$$

then $g$ cannot be degenerate w.r.t. $R$, and the function $\psi$ becomes strictly (log-)convex. Most importantly, we can conclude that $\partial_\theta \psi(0) < 0$.

LEMMA 5.8. *With $\pi_1^*$ and $\pi_2^*$ denoting the marginals of $\pi^*$ it holds that $\pi_1^* \otimes \pi_Q(f) < \mu^*$ as well as $\pi_P \otimes \pi_2^*(f) < \mu^*$.*

PROOF. We consider $(X_n, Y_n)_{n \geq 1}$ under the tilted measure $\mathbb{P}_{\pi^*}^*$ and an independent stationary Markov chain $(Z_n)_{n \geq 1}$ with transition probabilities $Q$. Then

$$(X_n, Y_n, Z_n)_{n \geq 1}$$

is a Markov chain on $E^3$, and we define the function $\tilde{f}$ on $E^3$ by

$$\tilde{f}(x, y, z) = f(x, z) - f(x, y).$$

The Markov chain has transition probabilities given by

$$R^*_{(x,y),(x',y')} Q_{z,z'} = \frac{r^*(x', y')}{r^*(x, y)} \exp(\theta^* f(x', y')) P_{x,x'} Q_{y,y'} Q_{z,z'},$$

with invariant distribution $\pi^* \otimes \pi_Q$. We also introduce the $\tilde{\Phi}(\theta)$ matrix

$$\tilde{\Phi}(\theta)_{(x,y,z),(x',y',z')} = \exp(\theta(f(x', z') - f(x', y'))) R^*_{(x,y),(x',y')} Q_{z,z'}.$$

With $\tilde{\varphi}(\theta)$ the spectral radius of $\tilde{\Phi}(\theta)$ we have that $\tilde{\varphi}(0) = \tilde{\varphi}(\theta^*) = 1$ since $\tilde{\Phi}(0)$ is stochastic and $\tilde{\Phi}(\theta^*)$ has a right eigenvector with eigenvalue 1 having entries $r^*(x,z)/r^*(x,y)$. Moreover, (8) provides the necessary cycle to show that (46) holds, and since

$$\partial_\theta \tilde{\varphi}(0) = \pi^* \otimes \pi_Q(\tilde{f}) = \pi_1^* \otimes \pi_Q(f) - \pi^*(f) = \pi_1^* \otimes \pi_Q(f) - \mu^*$$

by (9), it follows that $\pi_1^* \otimes \pi_Q(f) < \mu^*$. The second inequality follows similarly. □

LEMMA 5.9. *The sequence $(\mu_T^*)_{T \geq 1}$ is convergent, and*

$$\mu_\infty^* := \lim_{T \to \infty} \mu_T^* < \mu^*.$$

PROOF. We first observe that

$$\mu_T^* = \mathbb{E}_{\pi^*}^*(f(X_1, Y_{1+T})) \to \pi_1^* \otimes \pi_2^*(f)$$

for $T \to \infty$, where $\pi_1^*$ and $\pi_2^*$ are the marginals of $\pi^*$.

We consider $(X_n, Y_n)_{n \geq 1}$ under the tilted measure $\mathbb{P}_{\pi^*}^*$ and let $(W_n, Z_n)_{n \geq 1}$ be an *independent* copy with the same distribution. Then

$$(X_n, Y_n, W_n, Z_n)_{n \geq 1}$$

is a Markov chain on $E^4$ with transition probabilities $R^*_{(x,y),(x',y')} R^*_{(w,z),(w',z')}$ and invariant distribution $\pi^* \otimes \pi^*$. We define the function $f_\infty$ on $E^4$ by

$$f_\infty(x, y, w, z) = f(x, z) + f(w, y) - f(x, y) - f(w, z).$$

Introduce the corresponding $\tilde{\Phi}_\infty(\theta)$ matrix by

$$\tilde{\Phi}_\infty(\theta)_{(x,y,w,z),(x',y',w',z')} = \exp(\theta f_\infty(x',y',w',z'))R^*_{(x,y),(x',y')}R^*_{(w,z),(w',z')}$$

and its spectral radius $\tilde{\varphi}_\infty(\theta)$. By arguments analogous to those in Lemma 5.8 we conclude that $\tilde{\varphi}_\infty(0) = \tilde{\varphi}_\infty(\theta^*) = 1$, and that $\partial_\theta \tilde{\varphi}_\infty(0) = 2\mu^*_\infty - 2\mu^* < 0$. Hence $\mu^*_\infty < \mu^*$. Regarding $\partial_\theta \tilde{\varphi}_\infty(0) < 0$, we can again use (8) to verify that (46) holds. $\square$

It is interesting and very useful that the inequality in Lemma 5.9 holds not only in the limit but in fact for all $T$.

LEMMA 5.10.  *For all $T \geq 1$ it holds that*

(47) $$\mu^*_T < \mu^*.$$

PROOF.  With $S^T_n = \sum_{k=1}^n f(X_k, Y_{k+T})$ and $S_n = \sum_{k=1}^n f(X_k, Y_k)$ we observe that $S_n \overset{\mathcal{D}}{=} S^T_n$ under $\mathbb{P} = \mathbb{P}_\pi$, since the $X$- and $Y$-sequences are independent, stationary Markov chains. By (35) this implies that

(48) $$\frac{1}{n} \log \mathbb{E}(\exp(\theta S^T_n)) \to \log \varphi(\theta)$$

for $n \to \infty$.

Consider first the case $T = 1$ and the Markov chain

$$(X_n, X_{n+1}, Y_n, Y_{n+1})_{n \geq 1},$$

which under the tilted measure has transition probabilities

$$R^*_{(x,w,y,z),(x',w',y',z')} = \frac{r^*(w',z')}{r^*(w,z)} \exp(\theta^* f(w',z'))P_{w,w'}Q_{z,z'}\delta_{w,x'}\delta_{z,y'}.$$

Introduce the matrix

$$\tilde{\Phi}_1(\theta)_{(x,w,y,z),(x',w',y',z')} = \exp(\theta(f(x',z') - f(w',z'))) R^*_{(x,w,y,z),(x',w',y',z')}$$

and its spectral radius $\tilde{\varphi}_1(\theta)$. Clearly, $\tilde{\varphi}_1(0) = 1$ and we observe that

$$\tilde{\Phi}_1(\theta^*)_{(x,w,y,z),(x',w',y',z')} = \frac{r^*(w',z')}{r^*(w,z)} \exp(\theta^* f(x',z')) P_{w,w'} Q_{z,z'} \delta_{w,x'} \delta_{z,y'}.$$

The matrix $\tilde{\Phi}_1(\theta^*)$ has the same spectrum if we remove the eigenvector fraction, hence (35) together with (48) imply that

$$\log \tilde{\varphi}_1(\theta^*) = \lim_{n \to \infty} \frac{1}{n} \log \mathbb{E}(\exp(\theta^* S^1_n)) = \log \varphi(\theta^*) = 0,$$

thus $\tilde{\varphi}_1(\theta^*) = 1$.

Furthermore, by (9) $\partial_\theta \tilde{\varphi}_1(0) = \mu_1^* - \mu^*$. Using (8) (for $T = 1$) together with (46) we find that $\partial_\theta \tilde{\varphi}_1(0) < 0$, hence

$$\mu_1^* < \mu^*.$$

A similar argument for $T \geq 2$ is possible by introducing the Markov chain

$$(X_n, \ldots, X_{n+T}, Y_n, \ldots, Y_{n+T})_{n \geq 1}$$

and a function $f_T$ given by

$$f_T(x_0, \ldots, x_T, y_0, \ldots, y_T) = f(x_0, y_T) - f(x_T, y_T).$$

The spectral radius $\tilde{\varphi}_T(\theta)$ of the corresponding matrix $\tilde{\Phi}_T(\theta)$ fulfills that $\tilde{\varphi}_T(0) = \tilde{\varphi}_T(\theta^*) = 1$ and that $\partial_\theta \tilde{\varphi}_T(0) = \mu_T^* - \mu^* < 0$, using (8) to show that (46) holds. Thus $\mu_T^* < \mu^*$. □

5.6. *Variables shared in both sequences.* We define for $i, j, m, T \geq 1$ with $i \leq j$

$$S_1 = \sum_{k=1}^{i} f(X_k, Y_k), \qquad S_2 = \sum_{k=i+1}^{j} f(X_k, Y_k),$$

$$\widetilde{S}_2 = \sum_{k=i+1}^{j} f(X_k, Y_{k+T}), \qquad S_3 = \sum_{k=j+1}^{i+m} f(X_k, Y_{k+T}).$$

LEMMA 5.11. *There exist an $\varepsilon > 0$ and some $K$ (both independent of $T$) such that*

(49) $$\mathbb{P}(S_1 + S_2 > t, \widetilde{S}_2 + S_3 > t) \leq K \exp(-\theta^*(1 + \varepsilon)t)$$

*for $t \geq 0$.*

PROOF. Assume first that the number of variables $j - i$ in the overlapping part is small, less than $t(4\|f\|_\infty)^{-1}$, say, in which case we obtain the estimate

$$\mathbb{P}(S_1 + S_2 > t, \widetilde{S}_2 + S_3 > t) \leq \mathbb{P}(S_1 > 3t/4, S_3 > 3t/4)$$
$$\leq \rho \mathbb{P}(S_1 > 3t/4) \mathbb{P}(S_3 > 3t/4)$$
$$\leq K \exp(-3\theta^* t/2),$$

using Lemma 5.1 for the second inequality and then a standard exponential change of measure argument; see (36). This implies (49) with $\varepsilon = 1/2$.

If instead $j - i \geq t(4\|f\|_\infty)^{-1}$ we observe that

(50)
$$\mathbb{P}(S_1 + S_2 > t, \widetilde{S}_2 + S_3 > t)$$
$$\leq \mathbb{P}(S_1 + S_2 > t, \widetilde{S}_2 \geq S_2) + \mathbb{P}(\widetilde{S}_2 + S_3 > t, S_2 \geq \widetilde{S}_2).$$

With $L_j^* = r(X_j, Y_j)/r(X_0, Y_0) \exp(\theta^*(S_1 + S_2))$ we obtain

$$\mathbb{P}_\pi(S_1 + S_2 > t, \widetilde{S}_2 \geq S_2) = \mathbb{P}_\pi\left(\frac{L_j^*}{L_j^*}; S_1 + S_2 > t, \widetilde{S}_2 \geq S_2\right)$$

$$\leq \gamma \exp(-\theta^* t) \mathbb{P}_{\pi,j}^*(\widetilde{S}_2 \geq S_2),$$

where $\mathbb{P}_{\pi,j}^*$ denotes the tilted measure *up to index $j$*. Using Lemma 5.1, we can, at the expense of a factor $\rho$, assume that the sequence $(X_n, Y_n)_{n=i}^j$ is a stationary Markov chain under the tilted measure and that $(Y_n)_{n \geq j+1}$ is independent and stationary under the original measure. Under this assumption it follows that the mean of $\widetilde{S}_2 - S_2$ equals $(j - T - i)\mu_T^* + T\pi_1^* \otimes \pi_Q(f) - (j - i)\mu^*$. Using Lemmas 5.8, 5.9 and 5.10 we can find a $\zeta > 0$, independent of $T$, such that

$$(j - T - i)\mu_T^* + T\pi_1^* \otimes \pi_Q(f) - (j - i)\mu^* < -(j - i)\zeta.$$

Hence Theorem 5.7 gives that

$$\mathbb{P}_{\pi,j}^*(\widetilde{S}_2 \geq S_2) \leq \rho \exp\left(-\frac{\zeta^2(j-i)}{2\eta^2}\right) \leq \rho \exp\left(-\frac{\zeta^2 t}{8\|f\|_\infty \eta^2}\right)$$

or, with $\varepsilon = \zeta^2(\theta^* 8\|f\|_\infty \eta^2)^{-1}$,

$$\mathbb{P}(S_1 + S_2 \geq t, \widetilde{S}_2 \geq S_2) \leq \rho\gamma \exp(-\theta^*(1 + \varepsilon)t).$$

Of course, a similar argument takes care of the second term in (50) and (49) follows. □

5.7. *Useful mixing inequalities.* When the aligned sequences are i.i.d. the sets $B_a$ entering Theorem 4.1 are usually chosen such that $V_a$ and $\mathcal{F}_a$ are independent, in which case $\mathbb{E}|\mathbb{E}(V_a|\mathcal{F}_a) - \mathbb{E}(V_a)| = 0$ and the term $\beta_{4,n}$ in Theorem 4.1 vanish. In the framework of Markov chains we need to control $\beta_{4,n}$ by using exponential $\beta$-mixing of stationary, finite state-space Markov chains. To this end we need a few results on how to translate knowledge about the $\beta$-mixing coefficients into useful bounds on $\mathbb{E}|\mathbb{E}(V_a|\mathcal{F}_a) - \mathbb{E}(V_a)|$.

For two $\sigma$-algebras $\mathcal{F}$ and $\mathcal{G}$ the $\alpha$-mixing measure of dependence is

$$\alpha(\mathcal{F}, \mathcal{G}) = \sup_{A \in \mathcal{F}, B \in \mathcal{G}} |\mathbb{P}(A \cap B) - \mathbb{P}(A)\mathbb{P}(B)|.$$

The following lemma relates $\alpha$-mixing measures to mean values of the desired form.

LEMMA 5.12. *Let $\mathcal{F}$ and $\mathcal{G}$ be $\sigma$-algebras and let $A \in \mathcal{G}$. With $\eta = 1(A)$*

(51) $$\mathbb{E}|\mathbb{E}(\eta|\mathcal{F}) - \mathbb{E}(\eta)| \leq 2\alpha(\mathcal{F}, \mathcal{G}).$$

PROOF. With $B = (\mathbb{E}(\eta|\mathcal{F}) \geq \mathbb{E}(\eta)) \in \mathcal{F}$ and $\xi = 1(B)$ we see that

$$\mathbb{E}|\mathbb{E}(\eta|\mathcal{F}) - \mathbb{E}(\eta)| = \mathbb{E}(\xi(\mathbb{E}(\eta|\mathcal{F}) - \mathbb{E}(\eta))) - \mathbb{E}((1-\xi)(\mathbb{E}(\eta|\mathcal{F}) - \mathbb{E}(\eta)))$$
$$= 2(\mathbb{E}(\xi\eta) - \mathbb{E}(\xi)\mathbb{E}(\eta))$$
$$= 2(\mathbb{P}(A \cap B) - \mathbb{P}(A)\mathbb{P}(B)) \leq 2\alpha(\mathcal{F}, \mathcal{G}). \quad \square$$

The $\beta$-mixing measure of dependence between the $\sigma$-algebras $\mathcal{F}$ and $\mathcal{G}$ is defined as

$$\beta(\mathcal{F}, \mathcal{G}) = \mathbb{E}\left(\sup_{A \in \mathcal{F}} |\mathbb{P}(A|\mathcal{G}) - \mathbb{P}(A)|\right).$$

For a stationary stochastic process $(Z_n)_{n \in \mathbb{Z}}$ and for a subset $I \subseteq \mathbb{Z}$ we define the $\sigma$-algebra $\mathcal{F}_I = \sigma(Z_n; n \in I)$. The $\beta$-mixing coefficient is defined as

$$(52) \quad \beta(n) = \beta(\mathcal{F}_{[n,\infty)}, \mathcal{F}_{(-\infty,0]}) = \mathbb{E}\left(\sup_{A \in \mathcal{F}_{[n,\infty)}} |\mathbb{P}(A|\mathcal{F}_{(-\infty,0]}) - \mathbb{P}(A)|\right),$$

for $n \geq 1$ and the process $(Z_n)_{n \in \mathbb{Z}}$ is called $\beta$-mixing if $\beta(n) \to 0$ for $n \to \infty$. For two subsets $I, J \subseteq \mathbb{Z}$, the distance, $d(I, J)$, between the sets is defined as

$$d(I, J) = \inf_{n \in I, m \in J} |n - m|.$$

If $I, J \subseteq \mathbb{Z}$, we write $I < J$ if $n < m$ for all $n \in I$ and $m \in J$.

LEMMA 5.13. *Assume that $I_1 < J < I_2$ are three subsets of $\mathbb{Z}$. With $I = I_1 \cup I_3$ it holds that*

$$\alpha(\mathcal{F}_I, \mathcal{F}_J) \leq 3\beta(d(I, J)).$$

This result is Theorem 3.1 in [20]. See also [9], Theorem 1.3.3 for a slightly more general result.

5.8. *Proof of the Poisson approximation.* We recall the notation from Remarks 3.6 and 3.7 where

$$\tau_{-}(1) = \inf\{n > 0 | S_n \leq 0\}$$

and we let $\varepsilon_\delta = \varepsilon_{(0,0),\delta}$ such that $S_\delta = \delta\varepsilon_\delta(f)$.

LEMMA 5.14. *There exist constants $K, c > 0$ such that for all $n \geq 1$*

$$(53) \quad \mathbb{P}(\tau_{-}(1) \geq n) \leq K\exp(-cn).$$

*Moreover, for any $\eta > 0$ there exist constants $K(\eta), c(\eta) > 0$ such that for all $\delta \geq 1$*

$$(54) \quad \mathbb{P}(S_\delta > t, d(\varepsilon_\delta, \tilde{\pi}) \geq \eta) \leq K(\eta)\exp(-\theta^* t - c(\eta)\delta).$$

PROOF. We first note that there exists $\theta > 0$ such that $\log \varphi(\theta) < 0$ due to $\log \varphi(0) = 0$ and $\partial_\theta \log \varphi(0) = \partial_\theta \varphi(0) = \mu < 0$. Choose such a $\theta > 0$ and let $c = -\log \varphi(\theta) > 0$—for optimality we may choose $\theta$ that minimizes $\log \varphi(\theta)$. Exponential change of measure gives

$$\mathbb{P}(\tau_-(1) \geq n) = \mathbb{E}^{\theta f}\left(\frac{1}{L^{\theta f}_{\tau_-(1)}}; n \leq \tau_-(1) < \infty\right)$$

$$\leq K_0 \mathbb{E}^{\theta f}(\exp(-c\tau_-(1) - \theta S_{\tau_-(1)}); n \leq \tau_-(1) < \infty)$$

$$\leq K_0 \exp(-cn) \mathbb{E}^{\theta f}(\exp(-\theta S_{\tau_-(1)}); n \leq \tau_-(1) < \infty)$$

$$\leq K \exp(-cn).$$

Here $K_0$ is the maximum of the eigenvector fractions and

$$K = K_0 \mathbb{E}^{\theta f}(\exp(-\theta S_{\tau_-(1)})),$$

which is finite because $0 \geq S_{\tau_-(1)} \geq \min_{x,y} f(x,y)$. This shows (53).

For the second inequality we find that

$$\mathbb{P}(S_\delta > t, d(\varepsilon_\delta, \hat{\pi}) \geq \eta) = \mathbb{E}^*\left(\frac{1}{L^*_\delta}; S_\delta > t, d(\varepsilon_\delta, \hat{\pi}) \geq \eta\right)$$

$$\leq K_0 \exp(-\theta^* t) \mathbb{P}^*(d(\varepsilon_\delta, \hat{\pi}) \geq \eta).$$

Large deviation theory for Markov chains gives that

$$\limsup_{\delta \to \infty} \frac{1}{\delta} \log \mathbb{P}^*(d(\varepsilon_\delta, \hat{\pi}) \geq \eta) \leq - \inf_{\nu : d(\nu, \hat{\pi}) \geq \eta} I^2(\nu),$$

where the infimum is taken over all shift-invariant probability measures on $E^2 \times E^2$. Here the *rate-function* $I^2$ is continuous and $I^2(\nu) > 0$ for all $\nu \neq \hat{\pi}$. We refer to Definition III.23, Theorem IV.3 and Lemma IV.5 in [8]. Consequently we can choose $K_0(\eta), c(\eta) > 0$, with $c(\eta) < \inf_{\nu : d(\nu, \hat{\pi}) \geq \eta} I^2(\nu)$, such that for all $\delta \geq 1$

$$\mathbb{P}^*(d(\varepsilon_\delta, \hat{\pi}) \geq \eta) \leq K_0(\eta) \exp(-c(\eta)\delta).$$

We conclude that (54) holds with $K(\eta) = K_0 K_0(\eta)$. □

LEMMA 5.15. *If we, for some $x \in \mathbb{R}$, let*

(55) $$t = t_n = \frac{\log K^* + \log n^2 + x}{\theta^*}$$

*and assume that $(l_n)_{n \geq 1}$ is a sequence of positive integers satisfying*

$$\lim_{n \to \infty} l_n^{-1} \log n = \lim_{n \to \infty} n^{-1} l_n = 0,$$

*then with $x_n = \theta^*(t_n - \lfloor t_n \rfloor) \in [0, \theta^*)$ it holds that*

$$\sum_{a \in I} \mathbb{E}(V_a(t_n, l_n, \eta)) \sim \mathbb{E}(C_n(t_n)) \sim \exp(-x + x_n)$$

*for all $\eta > 0$.*

PROOF. Introduce the probabilities

$$p(n, x, y) = \mathbb{P}_{x,y}\left(\max_{\delta:\delta \leq \tau_-(1)} S_\delta > t_n\right) = \mathbb{P}_{x,y}\left(\max_{\delta:\delta \leq \tau_-(1)} S_\delta > \lfloor t_n \rfloor\right)$$

and

$$\tilde{p}(n, x, y) = \mathbb{P}_{x,y}\left(\max_{\substack{\delta:\delta \leq \tau_-(1) \wedge l_n \\ \text{and } d(\varepsilon_\delta, \hat{\pi}) \leq \eta}} S_\delta > t_n\right).$$

Furthermore, for $a = (i, j) \in I$ let

$$q(a, x, y) = \mathbb{P}(T_a = 0, X_i = x, Y_j = y).$$

Using the Markov property we find that the conditional probability of the event $(\max_{\delta:\delta \leq \Delta(a)} \delta \varepsilon_{a,\delta}(f) > t_n)$, conditionally on $(T_a = 0, X_i = x, Y_j = y)$, is smaller than $p(n, x, y)$ because $\Delta(a)$ is restricted by the boundaries of the score matrix. Thus

(56) $$\mathbb{P}\left(T_a = 0, \max_{\delta:\delta \leq \Delta(a)} \delta \varepsilon_{a,\delta}(f) > t_n\right) \leq \sum_{x,y \in E} p(n, x, y) q(a, x, y),$$

which by (23) gives

$$\mathbb{E}(C_n(t_n)) \leq \sum_{x,y \in E} p(n, x, y) \sum_{a \in I} q(a, x, y).$$

With $\tilde{I} = \{(i, j) \in I | i, j \leq n - l_n\}$ we find for $a = (i, j) \in \tilde{I}$, by conditioning on the event $(T_a = 0, X_i = x, Y_j = y)$, that

$$\mathbb{P}(V_a = 1) = \sum_{x,y \in E} \tilde{p}(n, x, y) q(a, x, y)$$

and hence

$$\sum_{a \in \tilde{I}} \mathbb{E}(V_a) = \sum_{x,y \in E} \tilde{p}(n, x, y) \sum_{a \in \tilde{I}} q(a, x, y).$$

Since by construction $\sum_{a \in I} V_a \leq C_n(t_n)$ we get the following chain of inequalities:

$$\sum_{x,y \in E} \tilde{p}(n, x, y) \sum_{a \in \tilde{I}} q(a, x, y) \leq \sum_{a \in I} \mathbb{E}(V_a)$$

$$\leq \mathbb{E}(C_n(t_n)) \leq \sum_{x,y \in E} p(n, x, y) \sum_{a \in I} q(a, x, y).$$

We are done once we have shown that the lower and upper bounds both behave as $\exp(-x + x_n)$. To this end, first note that by (19)

$$p(n, x, y) \exp(\theta^* \lfloor t_n \rfloor) \to e(x, y),$$

for $n \to \infty$ and as a consequence of (18), essentially considering the score matrix one diagonal at the time, we find that

$$\frac{1}{n^2} \sum_{a \in I} q(a, x, y) \to \frac{\nu(x, y)}{\mu^-}$$

for $n \to \infty$. This gives that

$$\exp(x - x_n) \sum_{x, y \in E} p(n, x, y) \sum_{a \in I} q(a, x, y)$$

$$= \frac{1}{K^*} \sum_{x, y \in E} p(n, x, y) \exp(\theta^* \lfloor t_n \rfloor)\, n^{-2} \sum_{a \in I} q(a, x, y)$$

$$\to \frac{1}{K^* \mu^-} \sum_{x, y \in E} e(x, y) \nu(x, y) = 1$$

for $n \to \infty$.

Regarding the lower bound, we observe that

$$\tilde{p}(n, x, y) \leq p(n, x, y)$$
$$\leq \tilde{p}(n, x, y) + \mathbb{P}(\tau_-(1) \geq l_n) + \mathbb{P}(\exists \delta \leq l_n : S_\delta > t_n, d(\varepsilon_\delta, \hat{\pi}) > \eta).$$

Since $l_n^{-1} \log n \to 0$ for $n \to \infty$ we conclude from (53) that $\mathbb{P}(\tau_-(1) \geq l_n) = o(\exp(\theta^* t_n))$. For the last probability on the right-hand side above we first observe that if $S_\delta > t_n$, then necessarily $\delta \geq \|f\|_\infty^{-1} t_n$. Thus using (54) and that $n^{-1} l_n \to 0$ for $n \to \infty$ we see that

$$\mathbb{P}(\exists \delta \leq l_n : S_\delta > t_n, d(\varepsilon_\delta, \hat{\pi}) > \eta) \leq l_n \exp(-(\theta^* + c(\eta) \|f\|_\infty^{-1}) t_n)$$
$$= o(\exp(\theta^* t_n)).$$

Hence

$$\tilde{p}(n, x, y) \exp(\theta^* \lfloor t_n \rfloor) \to e(x, y),$$

for $n \to \infty$. Since $n^{-1} l_n \to 0$ we also have that

$$n^{-2} \sum_{a \in \tilde{I}} q(a, x, y) \to \frac{\nu(x, y)}{\mu^-}$$

for $n \to \infty$. By an argument similar to that above

$$\exp(x - x_n) \sum_{x, y \in E} \tilde{p}(n, x, y) \sum_{a \in \tilde{I}} q(a, x, y) \to 1$$

for $n \to \infty$, and this completes the proof. □

LEMMA 5.16. *With $(t_n)_{n\geq 1}$ and $(l_n)_{n\geq 1}$ chosen as in Lemma 5.15, assuming in addition that*

$$\lim_{n\to\infty} n^{-\varepsilon} l_n = 0$$

*for all $\varepsilon > 0$, then under the assumptions in Theorem 3.1, the conditions in Theorem 4.1 are fulfilled for some $\eta > 0$ with*

$$\lambda_n = \exp(-x + x_n),$$

*that is,*

$$\left\| \mathcal{D}\left( \sum_{a \in I} V_a(t_n, l_n, \eta) \right) - \mathrm{Poi}(\exp(-x + x_n)) \right\| \to 0.$$

PROOF. We define the neighborhood of strong dependence, $B_a$ for $a \in I$, as follows. Define for $a = (i,j) \in I$

$$B_a^1 = \{(k,m) \in I \,|\, |k-i| \leq 2l_n\}, \qquad B_a^2 = \{(k,m) \in I \,|\, |m-j| \leq 2l_n\},$$

and $B_a = B_a^1 \cup B_a^2$.

Note that (36) provides the bound $\mathbb{E}(V_a) \leq K \exp(-\theta^* t_n)$, and since $|I| = n^2$ and $|B_a| \leq 4nl_n$, then

$$\sum_{a \in I, b \in B_a} \mathbb{E}(V_a)\mathbb{E}(V_b) \leq K' l_n n^{-1} \to 0$$

for $n \to \infty$. This shows that (26) holds.

We prove that (27) is fulfilled by splitting the set $B_a$ into three disjoint sets and, depending on the set, give a bound of $\mathbb{E}(V_a V_b)$ for $b$ in each of these sets. For $a \in I$ let

$$B_a = C_a \cup D_a^1 \cup D_a^2$$

with $C_a$, $D_a^1$ and $D_a^2$ being the disjoint sets

$$C_a = B_a^1 \cap B_a^2, \qquad D_a^1 = B_a^1 \setminus C_a, \qquad D_a^2 = B_a^2 \setminus C_a.$$

Consider the case $b \in C_a$ and $b \neq a$. Using (32) together with Lemma 5.11 we can find an $\varepsilon > 0$ such that

$$\mathbb{E}(V_a V_b) \leq l_n^2 K \exp(-\theta^*(1+\varepsilon)t_n).$$

Hence, observing that $\sum_{a \in I} |C_a| \leq 16 l_n^2 n^2$,

$$\sum_{a \in I, b \in C_a, b \neq a} \mathbb{E}(V_a V_b) \leq K' l_n^4 n^{-2\varepsilon} \to 0$$

for $n \to \infty$.

For $b \in D_a^1$ use (32) together with Corollary 5.4, which applies due to (12), to find $\eta, \varepsilon, K > 0$ such that
$$\mathbb{E}(V_a V_b) \leq K l_n^2 \exp(-(3/2 + \varepsilon)\theta^* t_n).$$
Since $|D_a^1| \leq 4 l_n n$ we conclude that
$$\sum_{a \in I, b \in D_a^1} \mathbb{E}(V_a V_b) \leq K' l_n^3 n^{-3\varepsilon} \to 0$$
for $n \to \infty$. The same bound is obtainable for $b \in D_a^2$ (cf. the comment after Corollary 5.4), and all in all we conclude that (27) is fulfilled.

The two-dimensional process $(X_n, Y_n)_{n \geq 1}$ is a stationary, irreducible Markov chain on a finite state space, hence we can extend it to a doubly infinite, stationary process $(X_n, Y_n)_{n \in \mathbb{Z}}$, which is exponentially $\beta$-mixing. The $\beta$-mixing coefficients therefore satisfy
$$\beta(n) \leq K \exp(-\gamma n)$$
for some constants $K, \gamma > 0$. For $a = (i, i) \in I$ we define $I_1 = (-\infty, i - l_n]$, $I_2 = [i+1, i+l_n]$, and $I_3 = [i + 2l_n + 1, \infty)$, for which $d(I_1 \cup I_3, I_2) = l_n + 1$. With $I = I_1 \cup I_3$ and $J = I_2$, then clearly $\mathcal{F}_a \subseteq \mathcal{F}_I = \sigma(X_n, Y_n | n \in I_1 \cup I_3)$ and $V_a$ is measurable w.r.t. $\mathcal{F}_J = \sigma(X_n, Y_n | n \in I_2)$. By Lemmas 5.12 and 5.13 it follows that
$$\mathbb{E}|\mathbb{E}(V_a | \mathcal{F}_a) - \mathbb{E}(V_a)| \leq 2\alpha(\mathcal{F}_I, \mathcal{F}_J) \leq 6\beta(l_n + 1) \leq K' \exp(-\gamma l_n).$$
For nondiagonal $a = (i, j) \in I$ we can shift the $X$-process by stationarity to reduce the problem to the previous one and thus to obtain the same bound. This bound implies that
$$\sum_{a \in I} \mathbb{E}|\mathbb{E}(V_a | \mathcal{F}_a) - \mathbb{E}(V_a)| \leq K' n^2 \exp(-\gamma l_n) \to 0$$
for $n \to \infty$. This shows that (28) holds, and combining the bounds obtained in this proof with Lemma 5.15 we see that Theorem 4.1 gives the result. □

REMARK 5.17. We have a little flexibility left in the choice of $(l_n)_{n \geq 1}$. It does not matter how we choose this sequence precisely, as it is only an intermediate, technical necessity for the proof. We just need to make sure that a sequence can be chosen with the desired properties—and this is indeed the case.

FINISHING THE PROOF OF THEOREM 3.1. Having proved Lemma 5.16 we only need to verify (24) according to Corollary 4.2. To this end we note that by construction $\sum_{a \in I} V_a \leq C_n(t_n)$, hence
$$\mathbb{P}\left(\sum_{a \in I} V_a \neq C_n(t_n)\right) \leq \mathbb{E}(C_n(t_n)) - \sum_{a \in I} \mathbb{E} V_a \to 0$$
for $n \to \infty$ by Lemma 5.15. □

LOCAL ALIGNMENT OF MARKOV CHAINS 35

**6. Concluding remarks.** As mentioned in the Introduction, the result is a generalization of that obtained by [6] for aligning independent i.i.d. sequences. The overall strategy of constructing a counting variable that approximates $C_n(t)$, and whose asymptotic behavior can be derived from [2], Theorem 1, is identical to the strategy employed in [6], though we have chosen an approximation whose relation to $C_n(t)$ seems more obvious. We have also chosen to use exponential change of measure arguments to obtain most of the needed inequalities, whereas Dembo, Karlin and Zeitouni [6] rely more on combinatorial and large deviation inequalities.

One major challenge was to find a appropriate generalization of condition (E′) in [6] for the i.i.d. case. First the condition given by (22) was obtained directly, but this condition is not able to completely retain the i.i.d. case. Fortunately the referees insisted that another attempt should be made to obtain the correct generalization of (E′). As it turned out, it is essential in the construction of the variables $V_a$ to require that the pair-empirical measure is close to $\hat{\pi}$. Although this does not affect the asymptotic behavior of the expectations $\mathbb{E}(V_a)$, it does provide bounds, as a result of Lemma 5.4, on the expectations $\mathbb{E}(V_a V_b)$ that seem unobtainable otherwise.

Another major challenge was the generalization of the part of the proof of Lemma 2 in [6] called case (c), where a smart permutation argument relying on exchangeability of i.i.d. variables was used. The solution presented here, which works for Markov chains, is an application of the Azuma–Hoeffding inequality for martingales as described in Sections 5.4 and 5.5.

Besides this an extra argument based on mixing inequalities was needed in order to take care of the $\beta_{4,n}$-term, which was not present in the i.i.d. case.

**Acknowledgment.** Thanks are due to the referees for useful comments and for encouraging me to find a natural extension in the Markov setup of condition (E′) in [6].## REFERENCES

[1] ALSMEYER, G. (2000). The ladder variables of a Markov random walk. *Probab. Math. Statist.* **20** 151–168. MR1785244
[2] ARRATIA, R., GOLDSTEIN, L. and GORDON, L. (1989). Two moments suffice for Poisson approximations: The Chen–Stein method. *Ann. Probab.* **17** 9–25. MR0972770
[3] ASMUSSEN, S. (2003). *Applied Probability and Queues*, 2nd ed. Springer, New York. MR1978607
[4] BUNDSCHUH, R. (2000). An analytic approach to significance assessment in local sequence alignment with gaps. In *RECOMB 00*: *Proceedings of the Fourth Annual International Conference on Computational Molecular Biology* 86–95. ACM Press, New York.
[5] ÇINLAR, E. (1975). *Introduction to Stochastic Processes.* Prentice–Hall Inc., Englewood Cliffs, NJ. MR0380912

[6] DEMBO, A., KARLIN, S. and ZEITOUNI, O. (1994). Limit distribution of maximal non-aligned two-sequence segmental score. *Ann. Probab.* **22** 2022–2039. MR1331214
[7] DEMBO, A. and ZEITOUNI, O. (1998). *Large Deviations Techniques and Applications*, 2nd ed. Springer, New York. MR1619036
[8] DEN HOLLANDER, F. (2000). *Large Deviations*. Amer. Math. Soc., Providence, RI. MR1739680
[9] DOUKHAN, P. (1994). *Mixing*. Springer, New York. MR1312160
[10] GROSSMANN, S. and YAKIR, B. (2004). Large deviations for global maxima of independent superadditive processes with negative drift and an application to optimal sequence alignments. *Bernoulli* **10** 829–845. MR2093612
[11] KARLIN, S. and DEMBO, A. (1992). Limit distributions of maximal segmental score among Markov-dependent partial sums. *Adv. in Appl. Probab.* **24** 113–140. MR1146522
[12] KINGMAN, J. F. C. (1961). A convexity property of positive matrices. *Quart. J. Math. Oxford Ser. (2)* **12** 283–284. MR0138632
[13] LEDOUX, M. and TALAGRAND, M. (1991). *Probability in Banach Spaces*. Springer, Berlin. MR1102015
[14] MERCIER, S. (2001). Exact distribution for the local score of one i.i.d. random sequence. *J. Comput. Biol.* **8** 373–380. Available at http://www.liebertonline.com/doi/abs/10.1089/106652701752236197.
[15] MERCIER, S. and HASSENFORDER, C. (2003). Distribution exacte du score local, cas markovien. *C. R. Math. Acad. Sci. Paris* **336** 863–868. MR1990029
[16] O'CINNEIDE, C. (2000). Markov additive processes and Perron–Frobenius eigenvalue inequalities. *Ann. Probab.* **28** 1230–1258. MR1797311
[17] SIEGMUND, D. and YAKIR, B. (2000). Approximate $p$-values for local sequence alignments. *Ann. Statist.* **28** 657–680. MR1792782
[18] SIEGMUND, D. and YAKIR, B. (2003). Correction: "Approximate $p$-values for local sequence alignments" [*Ann. Statist.* **28** (2000) 657–680 MR2002a:62140]. *Ann. Statist.* **31** 1027–1031. MR1792782
[19] STEELE, J. M. (1997). *Probability Theory and Combinatorial Optimization*. SIAM, Philadelphia. MR1422018
[20] TAKAHATA, H. (1981). $L_\infty$-bound for asymptotic normality of weakly dependent summands using Stein's result. *Ann. Probab.* **9** 676–683. MR0624695
[21] WATERMAN, M. S., ed. (1995). *Introduction to Computational Biology*. Chapman and Hall/CRC, Boca Raton, FL.

DEPARTMENT OF APPLIED MATHEMATICS
AND STATISTICS
UNIVERSITY OF COPENHAGEN
UNIVERSITETSPARKEN 5
DK-2100 COPENHAGEN
DENMARK
E-MAIL: richard@math.ku.dk